\documentclass[10pt,a4paper]{amsart}
\usepackage[english]{babel}
\usepackage[utf8x]{inputenc}
\usepackage{ucs}
\usepackage{tikz}
\usepackage{amsmath}
\usepackage{amsthm}
\usepackage{amsfonts}
\usepackage{amssymb}
\usepackage{hyperref}
\usepackage{dsfont}
\usepackage{mathrsfs}
\usepackage{eqnarray}
\write18{}
\usepackage{xcolor}
\usepackage{graphicx}
\usepackage{pdfpages} 
\usepackage[nothm]{thmbox}
\usepackage{esint}

\hypersetup{
    colorlinks=true,
    linkcolor=red,
    filecolor=magenta,      
    urlcolor=cyan,
}

\numberwithin{equation}{section}

\newcommand{\HH}{\mathbb{H}}
\newcommand{\ZZ}{\mathbb{Z}}
\newcommand{\RR}{\mathbb{R}}
\newcommand{\CC}{\mathbb{C}}
\newcommand{\NN}{\mathbb{N}}
\newcommand{\inte}[1]{\underset{#1}{\int}}
\newcommand{\finte}[3]{\underset{#1}{\fint} \ {#2} \ {#3}}
\newcommand{\supp}[1]{\mathrm{supp}(#1)}

\newcommand{\limiinf}[1]{\underset{#1}{\liminf}}

\newcommand{\somme}[1]{\underset{#1}{\sum}}
\newcommand{\limi}[1]{\underset{#1}{\lim}}
\newcommand{\tends}[1]{\underset{#1}{\longrightarrow}}
\newcommand{\equivaut}[1]{\underset{#1}{\sim}}
\newcommand{\fonction}[5]{\begin{array}{r r c l}
					#1 \hspace{2mm} : & #2 & \to & #3 \\
					& #4 & \mapsto & #5 \\
			   \end{array}}
\newcommand{\fonctionbis}[4]{\begin{array}{r c l}
					   #1 & \to & #2 \\
					 #3 & \mapsto & #4 \\
			   \end{array}}	  
\newcommand{\quotient}[2]{{\raisebox{.2em}{$#1$}\left/\raisebox{-.2em}{$#2$}\right.}}

\DeclareMathOperator{\vol}{vol}

\DeclareMathOperator{\Id}{Id}

\theoremstyle{definition}

\newtheorem{definition}[equation]{Definition}

\theoremstyle{theorem}

\newtheorem{theoreme}[equation]{Theorem}
\newtheorem{proposition}[equation]{Proposition}
\newtheorem{corollary}[equation]{Corollary}

\newtheorem{lemma}[equation]{Lemma}

\title{Orbital functions and heat kernels of Kleinian groups}

\author{Adrien Boulanger}

\address{Aix-Marseille Universit\'{e}, CNRS, Centrale Marseille, I2M, UMR 7373, 13453 Marseille, France}

\thanks{The author was partially founded by the ERC n°647133 'IChaos'.}

\begin{document}

\maketitle

\begin{abstract}
We study orbital functions associated to Kleinian groups through the heat kernel approach developed in \cite{artmoiheatcounting1}.
\end{abstract}

\section{Introduction}
We denote by $\HH^d$ the hyperbolic space of dimension $d \ge 2$. We call Poincaré group any discrete torsion free subgroup of the group of isometries of $\HH^d$ preserving the orientation. Following the historical terminology introduced by Poincaré, a Poincaré group acting on $\HH^3$ is called Kleinian and one acting on $\HH^2$ is called Fuchsian. If $\Gamma$ is a Poincaré group, the quotient space $\quotient{\HH^d}{\Gamma}$ is a hyperbolic manifold that we denote by $M_{\Gamma}$. Given a point $\tilde{p} \in \HH^d$ we denote by $p$ its image in $M_{\Gamma}$. Conversely, we will denote by $\tilde{p} \in \HH^d$ any lift of a point $p \in M_{\Gamma}$. Given a pair of points $\tilde{x},\tilde{y} \in \HH^d$ and $\rho > 0$ we define the orbital function as
	\begin{align*}
		 N_{\Gamma}(\tilde{x},\tilde{y},\rho)  : & = \sharp \{ \Gamma \cdot \tilde{y} \cap B(\tilde{x}, \rho) \} \\
		 	 & = \sharp \{ \gamma \in \Gamma \ , \ d(\tilde{x}, \gamma \cdot \tilde{y}) \le \rho \} \ ,
	\end{align*}
where the symbol $\sharp$ stands for the cardinality of a set. The function $N_{\Gamma}$ is $\Gamma$-invariant with respect to both its spatial variables, and therefore descends to a well defined function on $M_{\Gamma} \times M_{\Gamma} \times \RR_+$. We also call the resulting function \textbf{the orbital function} and, by abuse of notation, we keep denoting it by $N_{\Gamma}(x,y,\rho)$. \\

The first work, to the author's best knowledge, dedicated to understand the behaviour when $\rho \to \infty$ of orbital functions associated to Poincaré  groups goes back to Delsarte \cite{artdelsarte} in the early 40's. His approach to counting problems was extensively developed since. It was considered again with Huber \cite{arthuber2} who was the first to obtain a precise asymptotic of orbital functions associated to co-compact Fuchsian groups. Namely, in this case he obtained
	$$ N_{\Gamma}(x,y,\rho) \equivaut{\rho \to \infty} \frac{\pi \ e^{\rho}}{\vol(M_{\Gamma})} \ ,$$
where the term $\pi \ e^{\rho}$ has to be understood as the volume of a hyperbolic disk of radius $\rho$. Shortly after, Selberg \cite{artselberg2,artselberg} extracted what was essential in order for this approach to work and laid down the general framework of his celebrated trace formula. This led to extend Huber's result to Poincaré groups acting co-compactly. In the 70's Patterson \cite{artpattersontermeerror}, relying on the same method, improved Huber's result in two directions for Fuchsian groups: he replaced the assumption that the action is co-compact by the one that the quotient manifold has finite hyperbolic volume and he gave a precise control of the error term. Using more tools from Analysis, Lax-Phillips \cite{articlelax} and Patterson \cite{articlepatteron2} were able to extend this result to the setting of fundamental groups of geometrically finite hyperbolic manifolds, a class of groups encompassing the ones acting co-compactly. The most general theorem of the theory is due to Roblin \cite{memoireroblin} who gave a simple asymptotic for the orbital function for a large class of actions (without estimating the error term though). In the constant curvature setting, his result specialises as one can replace the geometrically finite assumption with the assumption that $M_{\Gamma}$ carries a finite Bowen-Margulis-Sullivan measure. The latter assumption is weaker since Sullivan showed \cite{artsullivanmesurebms} that if $\Gamma$ is geometrically finite then $M_\Gamma$ has a (unique) finite Bowen-Margulis-Sullivan measure. Let us mention, to conclude this historical overview, that Roblin's method relies on Margulis' dynamical ideas \cite{artmargulismixing} and was precisely initiated by Patterson \cite{arttersondensiteconforme} and Sullivan \cite{artsullivandensiteconforme, artsullivanmesurebms}. We refer to \cite{livremartine} for a more detailed exposition of the rich history of the study of orbital functions associated to groups acting on CAT(-1) spaces. \\

This article aims at investigating the behaviour of orbital functions associated to Kleinian groups whose quotients do not carry any finite Bowen-Margulis measures. There are very few cases where even rough estimates of the orbital function are obtained in this case. The only result in this direction, in the setting of Poincaré groups, is due to Epstein \cite{articleepsteinspectraltheory1, articleepsteinspectraltheory2, articleepsteinspectraltheory3} and deals with Abelian covers of compact hyperbolic manifolds. \\

Epstein's method relies on tools similar to those used in this article, namely the use of the spectral theory of the manifold $M_{\Gamma}$. The method is limited to Abelian covers of finite volume manifolds. Indeed, the fact that the covering group is Abelian allows them to use the finite measure of the underlying compact manifold, encoding totally the Abelian cover in the homology of the latter. It seems therefore hard to generalise this approach to study Poincaré groups associated to non Abelian covers of finite volume manifolds, and even harder for hyperbolic manifolds which does not cover any finite volume hyperbolic manifold at all. \\

In this article we pursue the investigation of the relation that heat kernels entertain with orbital functions of Poincaré groups. This approach was initiated by the author in his PhD dissertation \cite{thesemoi} which eventually led to \cite{artmoiheatcounting1}. It was inspired by Margulis' dynamical idea to relate geodesic flows and counting problems and by Sullivan's idea to relate Brownian motions and geodesic flows (in constant negative curvature). Therefore, counting problems and heat kernels should be related somehow. This article quantifies this relation through the following theorem. \\

We denote by $\mathcal{V}_d(\rho)$ the volume of any hyperbolic ball of radius $\rho$ in $\HH^d$ and by $p_{\Gamma}$ the heat kernel of $M_{\Gamma}$. We refer to Section \ref{secheatoperator} for a precise definition of the latter.

\begin{theoreme}[main theorem]
	\label{maintheo}
	Let $\Gamma$ be a Kleinian group. Assume that there exist $x,y \in \mathbb{H}^3$, $\alpha \ge 0$ such that for $t$ large enough we have
		$$p_{\Gamma}(x,y,t) \ge t^{-\alpha} \ ,$$ 
then for any $x \in \HH^3$ we have
		$$ N_{\Gamma}(x,x,\rho) \equivaut{\rho \to \infty} p_{\Gamma} \left(x,x,\frac{\rho}{2} \right) \ \mathcal{V}_3(\rho) \ .$$
\end{theoreme}

Our main theorem does not require any finite measure assumption, note however that the lower bound required on the heat kernel implies that the Kleinian group has maximal critical exponent. The approach taken in this article has the benefit of giving a precise asymptotic of the orbital function provided that we restrain ourselves to the on-diagonal study. Compare to the weaker upper bound obtained in \cite[Theorem 1.2]{artmoiheatcounting1} which is valid everywhere and for any Poincaré groups. We will give an overview of the proof of Theorem \ref{maintheo} in Section \ref{secoverview}. The main benefit of our main theorem is to reduce orbital function estimates to heat kernel ones. Notice that, contrary to the orbital function, the heat kernel has the semigroup property which makes its analysis easier. \\

In the proof, we only use the dimension 3 to prove Proposition \ref{propheatkernellip}. The proof of this proposition relies on the fact that we have an explicit and simple formula for the heat kernel in this case. Others formulas of the same type exist in all dimensions and one can most probably adapt the proof proposed here using the same method. However, the calculations are messy and it would have technically burdened this article a lot. As the main corollary of the above theorem concerns Kleinian groups, the author decided against it. Moreover, it seems to the author that the following statement, from which Proposition \ref{propheatkernellip} follows should hold. This would give a more satisfying proof to Proposition \ref{propheatkernellip} as it is  interesting on its own in the author's opinion. \\

\begin{question}
Let $M$ be a complete Riemannian manifold with pinched sectional curvature. Does there exist, for any compact set $K \subset M$, a constant $C_K > 0$ such that for any $t \ge 1$ and any $x, y \in K$ 
	$$ | \nabla_y \ln (p_M(x,y,t)) | \le C_K \ ? $$
\end{question}

A historically important class of examples of groups that do not carry any finite Bowen-Margulis-Sullivan measure are geometrically infinite finitely generated Kleinian group, also referred as \textbf{degenerate Kleinian groups}. This class of Kleinian groups was extensively studied since they are the main objects of Thurston's hyperbolisation theorem \cite{livremcmumu,livreotal2001hyperbolization} and several important problems (all theorems now) as Marden's tameness conjecture \cite{arttamenesscalegarigabai} (proved independently by Agol, unpublished), Ahlfors' conjecture (which follows from Marden's conjecture and \cite{articlecanary}) or the ending lamination conjecture \cite{artendinglamination}. For more details, we recommend the surveys \cite{surveycanarymarden,surveyminsky}. The only studies of their orbital functions was toward the understanding of the critical exponents of the group, which is defined as
 $$  \delta_{\Gamma} := \limi{\rho \to \infty} \ \frac{ \ln \big( N_{\Gamma}(x,y,\rho)  \big)}{\rho} \ . $$
Note that it is not clear that the limit is well defined \textit{a priori}: it was shown to exist in \cite{artroblinfoncorb}. Bishop and Jones  \cite{articlebishopjones,artbishopjonescriticalgeometricfiniteness}, following Sullivan \cite{articlesullivandemicylindre}, were able to show that $\delta_{\Gamma} =2$ for geometrically tame Kleinian groups. As a corollary of \cite{articlecanary} and the tameness theorem, all degenerate Kleinian groups are geometrically tame, which settles the study of their critical exponents. However, the critical exponent is a much weaker invariant than a precise asymptotic of the orbital function and none of the methods developed in the above mentioned work apply. Note that Epstein's result, discussed above, encompasses the case of $\ZZ$-covers associated to compact hyperbolic 3-manifolds which fibre over the circle, which are special cases of degenerate Kleinian groups. Theorem \ref{maintheo} was actually designed to prove the

\begin{corollary}
	\label{maincor}
	Let $\Gamma$ be a degenerate Kleinian group such that $M_{\Gamma}$ has positive injectivity radius. Then the conclusion of Theorem \ref{maintheo} holds. Moreover, 
		\begin{itemize}
			\item either the limit set of $\Gamma$ is the entire sphere and there are positive constants $C_-, C_+, \rho_0$ such that for every $x \in M_{\Gamma}$ there is $\rho_0 > 0$ such that for $\rho \ge \rho_0$ we have
				$$  \frac{ C_- \ e^{ 2 \rho}}{\rho^{\frac{1}{2}}} \le N_{\Gamma}(x,x,\rho) \le \frac{ C_+ \ e^{2 \rho}}{\rho^{\frac{1}{2}}} \ ;$$
			\item or for any $x \in M_{\Gamma}$ there are two positive constants $C_-(x), C_+(x)$ such that for $\rho \ge 1$
				$$ \frac{ C_-(x) \ e^{ 2 \rho}}{\rho^{\frac{3}{2}}} \le N_{\Gamma}(x,x,\rho) \le  \frac{ C_+(x) \ e^{2 \rho}}{\rho^{\frac{3}{2}}} \ .$$
		\end{itemize}
\end{corollary}

The proof of the above corollary is an immediate application of Theorem \ref{maintheo} combined with the heat kernel estimates obtained in \cite{artmoiheatcounting1}. Indeed, under the same assumptions than in Corollary \ref{maincor} we showed \cite[Theorem 1.10]{artmoiheatcounting1} that 
\begin{itemize}
			\item either the limit set of $\Gamma$ is the entire sphere and there is two positive constants $C_-, C_+$ such that for every $x \in M_{\Gamma}$ and $t$ large enough we have 
				$$   \frac{C_-}{t^{\frac{1}{2}}} \le p_{\Gamma}(x,x,t) \le  \frac{C_+}{t^{\frac{1}{2}}} \ ;$$
			\item or for any $x \in M_{\Gamma}$ there is two positive constants $C_-(x), C_+(x)$ such that for $t \ge 1$ we have
				$$  \frac{ C_-(x)}{t^{\frac{3}{2}}} \le p_{\Gamma}(x,x,t) \le \ \frac{C_+(x)}{t^{\frac{3}{2}}} \ .$$
		\end{itemize}

The heat kernel can be interpreted as the transition kernels of the Brownian motion. This makes it more intuitive than the orbital function: we refer to \cite{artmoiheatcounting1} for a justification \textit{a priori} of the exponents $1/2$ and $3/2$. Note that we do not use the same terminology in \cite{artmoiheatcounting1}: the 'fully degenerate' case corresponds to the case where the limit set is the entire sphere and the 'mixed type' case to the other one. \\ 

The next natural step to investigate orbital functions of finitely generated Poincaré groups is to drop the assumption that $M_{\Gamma}$ has positive injectivity radius. Since McMullen conjecture on polynomial growth of degenerate ends was solved in \cite[volume growth theorem, page 5]{artendinglamination}, a positive answer to the following question implies that Theorem \ref{maintheo} holds for a broader class of degenerate Kleinian groups. \\

\begin{question}
	Given a Riemannian manifold $M$ with pinched sectional curvature carrying at least one infinite volume end of (strictly) polynomial growth, are there a point $x \in M_{\Gamma}$ and $\alpha > 0$ such that for $t$ large enough we have 	
		$ p_{M}(x,x,t) \ge t^{-\alpha} \ ?$ \\
\end{question}

\textbf{About the notation.} We will often change of perspective in this article, going back and forth from $\HH^d$ to $M_{\Gamma}$. As already mentioned, we will not make any difference, notation-wise, between kernels of operators acting on functions of $M_{\Gamma}$ or as $\Gamma \times \Gamma$-invariant kernels defined on $\HH^d \times \HH^d$. We shall make clear, when ambiguous, which space is the source by using the tildes to emphasis that we work on $\HH^d$. \\

\textbf{Acknowledgement.} I am very grateful to Gilles Courtois, my former PhD advisor, for his proof checking time as well for his many useful comments about the redaction of this article. I would also like to thank Sergiu Moroianu for interesting conversations around Selberg's trace formula and Alexander Grigor'yan for emails exchange related to heat kernels.

\section{First definitions and overview of the argument.}
\label{secoverview}
In order to prove our main Theorem \ref{maintheo} we start off introducing orbital operators, defined bellow. The operator approach is not the novelty of this article. First, we introduce the average operators, which are the underlying 'universal' operators behind the construction of the orbital ones. \\

We denote by $\mu_h$ the hyperbolic measure of any hyperbolic manifold and by $\fint_{E} f d \ \mu_h $ the mean value of the function $f$ over the set $E$ with respect to the measure $\mu_h$. We denote by $\pi : \HH^d \to M_{\Gamma}$ the universal Riemannian covering. \\

We call \textbf{average operator} the family of operators indexed by $\rho > 0$ acting on bounded functions on $\HH^3$ defined as 
		$$ \mathcal{O}^{\rho}(f)(\tilde{x}) := \finte{B(\tilde{x},\rho)}{f}{ d \mu_h} \ . $$

It is easy to see that if $f$ is $\Gamma$-invariant then  $\mathcal{O}^{\rho}(f)$ is also $\Gamma$-invariant. In particular, if $\tilde{f}$ is the lift of a function $f$ defined on $M_{\Gamma}$, the function $\mathcal{O}^{\rho}(\tilde{f})$ induces a well defined function on $M_{\Gamma}$ that we denote by  $\pi_*\mathcal{O}^{\rho}(\tilde{f})$. \\

Given a Riemannian manifold $(M,g)$, we denote by $B^+_0(M)$ the vector space of measurable non negative bounded and compactly supported functions. Note that we did not require the functions of $B^+_0(M)$ to be continuous. \\

We denote by $\supp{f}$ the support of a function $f$. Note that we have $B^+_0(M) \subset L^2(M)$, where $L^2(M)$ stands for the space of square integrable functions defined on $M$ with respect to the Riemannian measure. For a pair of functions $f,h \in L^2(M)$ we denote by
	$$ \left<f,h \right>_{L^2(M)} := \inte{M} \ f \ h \ d\mu_g $$
the usual scalar product on $L^2(M)$. \\

Note that given a function $f \in B^+_0(M_{\Gamma})$ and $\rho > 0$, the support of $\pi_*\mathcal{O}^{\rho}( \tilde{f})$ is contained in a $\rho$-neighbourhood of $\supp{f}$. In particular $\pi_*\mathcal{O}^{\rho}(f) \in B^+_0(M_{\Gamma})$.

\begin{definition}
	Given a Poincaré group $\Gamma$, we call \textbf{orbital operator}, that we denote by $(\mathcal{O}_{\Gamma}^\rho)_{\rho > 0}$, the family of operators acting on $B^+_0(M_{\Gamma})$ defined as 
	$$ \mathcal{O}_{\Gamma}^{\rho}(f) := \pi_* \mathcal{O}^{\rho}(\tilde{f}) \ . $$  
\end{definition} 

The proof of the following well known lemma is elementary and left to the reader. Recall that we denoted by $\mathcal{V}_d(\rho)$ the volume of a hyperbolic ball in $\HH^d$ of radius $\rho$.
\begin{lemma}
	\label{lemmaaverage}
	Let $\Gamma$ be a Poincaré group. The function 
	$$ \frac{N_{\Gamma}(\cdot, \cdot, \rho)}{\mathcal{V}_d(\rho)} $$
is the kernel of the orbital operator $\mathcal{O}_{\Gamma}^{\rho}$. Namely, for any $\rho > 0$ and any $f \in B^+_0(M_{\Gamma})$ we have 
	$$ \mathcal{O}_{\Gamma}^{\rho}(f)(x) = \inte{M_{\Gamma}} \ \frac{N_{\Gamma}(x,y,\rho)}{\mathcal{V}_d(\rho)}  \ f(y) \ d \mu_h(y) \ . $$  
\end{lemma}

In particular, the orbital function being symmetric in $x$ and $y$, the orbital operator is also symmetric with respect to the $L^2$-scalar product (in restriction to $B^+_0(M_{\Gamma}))$. \\

We will prove our main theorem thanks to the following proposition. We refer to Section \ref{secheatoperator} for all the material and definitions about the heat kernel.

\begin{proposition}[operators comparison]
	\label{propmaintheo}
	Let $\Gamma$ be a Poincaré group such that there is $\alpha > 0$, $x,y \in \mathbb{H}^d$ such that for $t$ large enough we have
		$$p_{\Gamma}(x,y,t) \ge t^{-\alpha} \ .$$ 
Then, for any non zero $f \in B^+_0(M_{\Gamma})$, 
		$$ \left< \mathcal{O}_{\Gamma}^{\rho}(f) , f \right>_{L^2(M_{\Gamma})} \equivaut{\rho \to \infty} \left< e^{- \frac{\rho}{2} \Delta}(f) , f \right>_{L^2(M_{\Gamma})} \ .$$
\end{proposition}

Roughly, since we know the heat operator and the orbital operator to commute in the setting of rank one symmetric spaces, the idea behind the proof of the above proposition is that two commuting operators can be compared only looking at their spectrum. However, since $M_{\Gamma}$ is not compact, its spectral theory is not easy to deal with, which makes hard to readily implement such an approach. In order to fall back in a simpler spectral setting, we shall approach $M_{\Gamma}$ with large domains $(\Omega_n)_{n \in \NN}$. The Dirichlet heat problem for bounded domains also has a unique solution, see Section \ref{secheatoperator}. In particular, a bounded domain carries a well defined heat kernel as well. We will recover the information needed about the heat kernel of the entire manifold $M_{\Gamma}$ using a theorem of Dozdiuk (Theorem \ref{theododziuk} Section \ref{secheatoperator}) which relates it to those of the sequence $(\Omega_n)_{n \in \NN}$.  \\

It is well known that relatively compact domains carries an $L^2$-eigenbasis of smooth eigenfunctions of the Laplace operator (Theorem \ref{theospectral}). This will be crucial for us. The following proposition asserts that the orbital operator and the heat one 'commute' for points far away of the boundary of a given domain. It is a variation of a well known fact. However, for comprehensiveness and since the author did not find any ready-to-use reference for domains, we give its (rather elementary) proof in the appendix of this article. Recall that the upper half plane model of the hyperbolic space $\mathbb{H}^d$ is given by $\RR^{2} \times \RR_+$ endowed with the metric 
	$$ \frac{dx^2_1 + ... + dx^2_{d-1} + dy^2}{y^{2}} \ . $$
\begin{proposition}[Delsarte's formula]
	\label{propselbergtransform} 
Let $\Gamma$ be a Poincaré group, $\Omega$ any open subset of $M_{\Gamma}$ and $\Phi$ a smooth function of $\Omega$ verifying $\Delta \Phi = \lambda \Phi$. Then for any $x \in \Omega$ such that $B(x,\rho) \subset \Omega$ we have
\begin{equation*}
	 \mathcal{O}_{\Gamma}^{\rho}(\Phi)(x) = \nu_{\rho}(\lambda) \ \Phi(x) \ ,
\end{equation*}

with 
\begin{equation}
	\label{eqformuleexplicit}
\nu_{\rho}(\lambda) = \mathcal{O}^{\rho}(y^s) = \finte{B(o,\rho)}{ y^s}{  \frac{dx_1...dx_{d-1} dy}{y^d}}  
\end{equation}
where $o := (0,..., 0, 1) \in \RR^{d-1} \times \RR_+^*$ and $s$ satisfying $s(d-1-s) = \lambda$.
\end{proposition}

 The function $\nu_{\rho}$ is usually referred to as a \textbf{Selberg's transform}. Note that if $f \in B^+_0(M_{\Gamma})$ then $ \supp{ \mathcal{O}^{\rho}_{\Gamma}(f)} \subset \Omega$ if $\Omega$ contains a $\rho$-neighbourhood of $\supp{f}$. We fix $\Omega$ with this property. The next step is to expand $\mathcal{O}^{\rho}_{\Gamma}(f)$ with respect to a $L^2$-eigenbasis of eigenfunctions $(\Phi_{\lambda})_{\lambda \in \mathrm{Sp}(\Omega)}$ associated to the Laplace operator of $\Omega$. Thanks to the above proposition, we will identify the coefficients associated to $\Phi_{\lambda}$ of this expansion to $\nu_{\rho}(\lambda)$. Informally, for any $f \in B^+_0(M_{\Gamma})$ and any domain $\Omega$ such that $\supp{\mathcal{O}^{\rho}_{\Gamma}(f)} \subset \Omega$ 
	$$ \mathcal{O}^{\rho}_{\Gamma}(f) = \somme{ \mathrm{\lambda \in Sp}(\Omega)} \nu_{\rho}(\lambda) \left< \Phi_{\lambda}, f \right>_{L^2(M)} \Phi_{\lambda} \ . $$
The Selberg's transform $\nu_{\rho}(\lambda)$ only depends on $\lambda$ and not on the underlying eigenfunction $\Phi_{\lambda}$. In particular it does not depend on $\Omega$ either (provided that it was taken large enough). This uniformity is crucial: it is ultimately what is allowing us to use the exhaustion of $M_{\Gamma}$ by relatively compact domains. Note that we do not prove anything about the spectral theory of the orbital operator itself. \\

It remains then to relate precisely the coefficients appearing in the expansion of $\mathcal{O}^{\rho}_{\Gamma}(f)$ (of the form $\nu_{\rho}(\lambda)$) to their expected values for the heat operator (of the form $e^{- \lambda \rho}$). This is what the following proposition aims at.  The proof of this proposition is rather technical (we split it in three steps) and will be carried out through Section \ref{secetudevpselberg}. Practically, this is a computation starting off Formula \eqref{eqformuleexplicit}. Heuristically, it is a spectrally quantified version of Sullivan's idea that Green functions associated to Brownian motion are related with critical exponents in the setting of rank 1 symmetric spaces and their quotients. See for example \cite{articlesullivanpositivity}. With this picture in mind, the term $1/(d-1)$ is natural and due to the fact that the escape rate of the Brownian motion is $d-1$ in $\HH^d$. 

\begin{proposition}[spectral proposition]
\label{propspectralprop}
For any $\beta >0$ there is a constant $C > 0$ such that for all $\epsilon > 0$ there is $\rho_0 > 0$ such that for all $\rho > \rho_0$ we have 
\begin{enumerate} 
	\item for all $ \lambda \le \frac{\beta \ln (\rho) }{\rho} $ 
			$$ | \nu_{\rho}(\lambda) - e^{-\lambda \frac{\rho}{d-1}}| \le \epsilon \ e^{ - \lambda \frac{\rho}{d-1}}  \ .$$ 
	\item for all $ \lambda \ge \frac{\beta \ln (\rho) }{\rho}$ 
		$$  |\nu_{\rho}(\lambda)|  \le C \ \rho^{-\frac{\beta}{d-1}} $$ 
\end{enumerate}
\end{proposition}

The first item quantifies that 'large' coefficients appearing in the expansion of $\mathcal{O}^{\rho}_{\Gamma}(f)$ are close to those associated to the heat operator. The second item, together with the lower bound assumption we made on the heat kernel, will imply that 'small' coefficients are negligible. \\

The operators comparison done, it will remain to relate them with their kernels. For the orbital operator part, we shall circumvent the analytical intricacies one may encounter with non smooth kernels by using Lemma \ref{lemmapproxorbital}, lemma which is inspired by an argument of Eskin and McMullen proposed in \cite{arteskinmac}. Note that Lemma \ref{lemmapproxorbital} together with the well known spectral Theorem \ref{theospectral} and the appendix of this article is enough to recover the asymptotic of orbital functions associated to groups acting by isometries and co-compactly on rank 1 symmetric spaces (even off diagonal). It simplifies the existing proofs in these cases but does not seem appropriate to get finer estimates. For the heat operator part, we rely on an explicit formula for the heat kernel of $\HH^3$. The 'local' regularity of $p_{\Gamma}$ will be precisely investigated though Proposition \ref{propheatkernellip}. 

\section{Heat operators and their kernels}

\label{secheatoperator}

\subsection{Heat kernels on Riemannian manifolds.} General references for this section are \cite{livrechavel} or \cite{bookgrigorheatkernel}. If $(M,g)$ is a Riemannian manifold, we will denote by $\mu_g$ the Riemannian measure and by $\Delta$ the Laplace operator with the convention so that $\Delta$ is a positive operator on $L^2(M)$.

\begin{theoreme}{\cite[Corollary 4.11, page 117]{bookgrigorheatkernel}}
\label{theocauchyprobleml2}
	Let $(M,g)$ be a complete Riemannian manifold. For any initial condition $u_0 \in L^2(M)$ the following Cauchy problem 
	\begin{equation}
		\label{eq Cauchy problem}
		 \left\{ 
			\begin{array}{l}
		 		\Delta u(x,t) + \partial_t u(x,t) = 0 \\
		 		u(\cdot,t) \tends{t \to 0} u_0 \hspace{0.5 cm} \text{ in \ }  \ L^2(M) \ . 
			\end{array} 
		 \right. 
	\end{equation}
	
has a unique solution $u : M \times \RR_+ \to \RR$ given by 
 $$u(t,x) = e^{-t \Delta}(u_0)(x,t) = \inte{M} \ p_M(x,y,t) \ f(y) \ d \mu_{g}(y) \ , $$
where $p_M(x,y,t)$ is, by definition, the \textbf{heat kernel} of $(M,g)$. 
\end{theoreme}

We will denote for short by $p_{\Gamma}(x,y,t)$ the heat kernel of the hyperbolic manifold $M_{\Gamma}.$ \\

There is a similar Cauchy problem for domains on Riemannian manifold. In order to ensure uniqueness one must require boundary conditions.

\begin{theoreme}{\cite[page 168]{ livrechavel}}
\label{theocauchyprobleml2}
	Let $\Omega$ be a bounded domain of a complete Riemannian manifold with smooth boundary and $u_0$ a continuous function on $\bar{\Omega}$ with $u_0 = 0$ on $\partial \Omega$. The following Cauchy problem 
	\begin{equation}
		\label{eq Cauchy problem}
		 \left\{ 
			\begin{array}{l}
		 		\Delta u(x,t) + \partial_t u(x,t) = 0 \\
		 		u(x,t) \tends{t \to 0} u_0(x) \text{ for every } \ x \in \bar{\Omega} \\
		 		u( \cdot ,t) = 0 \ \text{ on } \ \partial \Omega \ \text{ for } \  t > 0  \ .
			\end{array} 
		 \right. 
	\end{equation}
	
has a unique smooth solution $u(x,t) : \Omega \times \RR_+^* \to \RR$ given by 
 $$u(t,x) = e^{-t \Delta_{\Omega}}(u_0)(x,t) = \inte{M} \ p_M(x,y,t) \ u_0(y) \ d \mu_{g}(y) \ , $$
where $p_{\Omega}(x,y,t)$ is, by definition, the \textbf{heat kernel} of $\Omega$. 
\end{theoreme} 

Given an exhaustion by relatively compact open subsets $(\Omega_n)_{n \in \NN}$ with smooth boundaries of a manifold $M$, one can wonder how relate the heat kernel of the whole manifold $M$ with those associated to $(\Omega_n)_{n \in \NN}$. The following theorem gives a satisfying answer to this question.

\begin{theoreme}{\cite[Theorem 3.6]{artdodziuk}}
	\label{theododziuk}
	Let $M$ be a complete Riemannian manifold and $(\Omega_n)_{n \in \NN}$ an exhaustion of $M$ by open relatively compact subsets with smooth boundary. Then 
	$$p_{\Omega_n} \tends{n \to \infty} p_M$$
uniformly on every compact set of $\RR^*_+ \times M \times M$.
\end{theoreme}

The following well known theorem will provide us with a crucial tool of our approach, eigenfunctions of the Laplace operator. 

\begin{theoreme}{\cite[Theorem 10.13]{bookgrigorheatkernel}, \cite[page 169]{livrechavel}}
	\label{theospectral}
	Let $M$ be a compact Riemannian manifold (possibly with smooth boundary), then there exists a sequence $\mathrm{Sp}(M)$ of positive numbers going to infinity together with a family of smooth functions $(\Phi_{\lambda})_{\lambda \in Sp(M)}$ (vanishing on the boundary if any) in $L^2(M)$ of unit norm such that: 
	
	\begin{itemize}
		\item for any $\lambda \in \mathrm{Sp}(M)$ we have $\Delta \Phi_{\lambda} = \lambda \ \Phi_{\lambda}$;
		\item the family $(\Phi_{\lambda})_{\lambda \in Sp(M)}$ is a Hilbert basis of the space $L^2(M)$. Namely, for any $f \in L^2(M)$ we have the following convergence in $L^2(M)$
			$$ f =  \somme{\lambda \in \mathrm{Sp} (M)} \ \left<\Phi_{\lambda},f \right>_{L^2(M)} \Phi_{\lambda} \ .$$
	\end{itemize}
\end{theoreme}

\subsection{Local control of the heat kernel.} The two following statements address the regularity of the heat kernel and will be used to recover our main theorem from Proposition 'operators comparison'.

\begin{lemma}
	\label{lemmeheatdecreases}
	Let $M$ be a complete Riemannian manifold, then for every $x \in M$ the function 
	$$ \fonctionbis{\RR^*_+}{\RR_+}{t}{p(x,x,t)} $$
is non increasing.
\end{lemma}

\textbf{Proof.} Given a function $f \in L^2(M)$ we denote by $$\Phi_f(t) := \left< e^{-t\Delta}(f),f \right>_{L^2(M)} \ .$$ 
In order to deduce the conclusion of Lemma \ref{lemmeheatdecreases} we first show that for any $f$ as above we have
\begin{equation}
	\label{eqlemmeheatdec1}
		 \Phi_f'(t) \le 0 \ . 
\end{equation}

We compute the derivative 
	$$ \Phi_f'(t) = \left< \partial_t \left(e^{-t\Delta} (f)\right),f \right>_{L^2(M)} \ .$$ 

As a solution of the heat equation $e^{-t\Delta} (f)$ satisfies 
\begin{align*}
	\partial_t \left( e^{-t\Delta} (f) \right) & = - \Delta \left( e^{-t\Delta} (f) \right) \\  
		& = - e^{-\frac{t}{2} \Delta} \Delta \left( e^{-\frac{t}{2}\Delta} (f) \right) \ ,
\end{align*}

since the Laplace operator and the heat operator commute. Therefore, since the heat kernel is symmetric, we have
	\begin{align*}
		- \Phi_f'(t) & = \left< e^{-\frac{t}{2} \Delta} \Delta \left( e^{-\frac{t}{2}\Delta} (f) \right),f \right>_{L^2(M)} \\
		 & =  \left< \Delta \left( e^{-\frac{t}{2}\Delta} (f) \right), e^{-\frac{t}{2} \Delta}(f) \right>_{L^2(M)} \\
		 & = \left< \nabla \left( e^{-\frac{t}{2}\Delta} (f) \right), \nabla \left( e^{-\frac{t}{2} \Delta}(f) \right) \right>_{L^2(M)} \\
		 & = \big| \big| \nabla \left( e^{-\frac{t}{2}\Delta} (f) \right) \big| \big|_{L^2(M)} \ge 0 \ .
	\end{align*}

We conclude by localising Equation \eqref{eqlemmeheatdec1}. For any $x \in M_{\Gamma}$ and $\epsilon > 0$, we set
$$ \text{\large{X}}_{\epsilon} := \frac{\mathds{1}_{B(x,\epsilon)}}{\mathcal{V}_d(\epsilon)} \ . $$
though as an $L^2(M_{\Gamma})$-approximation of the Dirac mass at $x$. \\

Using \eqref{eqlemmeheatdec1} we get for any $\epsilon > 0$ and any $t \le s$
$$\Phi_{\text{{X}}_{\epsilon}}(t) \ge \Phi_{\text{{X}}_{\epsilon}}(s) \ .$$ 

Letting $\epsilon \to 0$ in the above inequality gives the conclusion using that heat kernels are continuous. \hfill $\blacksquare$ \\

The first item of the following statement will be needed in order to control the heat kernel in a neighbourhood of the diagonal only using the on-diagonal values. This is one of the two keys (with Lemma \ref{lemmapproxorbital}) which allows us to recover our main theorem from Proposition 'operator comparison'. The second item will be needed in order to improve the punctual polynomial lower bound assumed in our main Theorem to a uniform lower bound on a given compact sets. Practically, we will use it with the support of a function $f \in B^+_0(M)$ along the proof of Proposition 'operators comparison'.

\begin{proposition}
	\label{propheatkernellip}
	 Let $\Gamma$ be a Kleinian group fort which there exist $x,y \in M_{\Gamma}$ and $\alpha > 0$ such that for $t$ large enough 
			$$ p_{\Gamma}(x,y,t) \ge  t^{-\alpha} \ .$$ 
Then, for all $\epsilon > 0$ there is $\delta > 0$ such that
	\begin{enumerate} 
		\item \textbf{local version:} there is a time $t(x,y, \alpha) > 0$ such that for all $z \in B(y,\delta)$ and $t \ge t(x,y, \alpha)$ we have 
			$$(1 - \epsilon) \ p_{\Gamma}(x,y,t) \le p_{\Gamma}(x,z,t) \le (1 + \epsilon) \ p_{\Gamma}(x,y,t) \ ;$$
		\item \textbf{compact version}: for every compact set $K \in M_{\Gamma}$ there is a constant $C_K > 0$ and a time $t_K > 0$ such that for all $w,z \in K$ and all $t > t_K$ we have 
			$$ p_{\Gamma}(w,z,t) \ge C_K \  t^{-\alpha}\ . $$	
	\end{enumerate}
\end{proposition}

The next subsection is dedicated to the proof of the above proposition, which is rather technical. 

\subsection{Proof of Proposition \ref{propheatkernellip} } The main ingredient of the proof is the following theorem which gives an explicit formula for the heat kernel of $\HH^3$.

\begin{theoreme}[see for example \cite{artgrigoheatkernelhyper}]
\label{theogrigoformuleheatkernelhyp}
Let $x, y \in \HH^3$ and $t > 0$, then 
\begin{equation}
	\label{formuleniyauchaleur}
	 p_{\HH^3}(x,y,t) = \frac{1}{(4\pi t)^{\frac{3}{2}}} \frac{\rho}{\sinh(\rho)} e^{- t - \frac{\rho^2}{4t}} \ ,
\end{equation}
 where $\rho = d(x,y)$.
 \end{theoreme}
 
The above formula shows that the heat kernel only depends on the distance between $x$ and $y$. This is actually more general since it always holds for rank 1 symmetric spaces. We denote by $p_{3}(\rho,t)$  the heat kernel on $\HH^3$ viewed as a function of the distance. \\

We shall now relate the heat kernel of $\HH^3$ to the one of $M_{\Gamma}$ through the following well known lemma.

\begin{lemma}
	\label{lemheatkernelseries}
Let $\Gamma$ be a Kleinian group. Then for any $x, y \in M_{\Gamma}$ and any $t > 0$ we have 
	$$ p_{\Gamma}(x,y,t) = \somme{\gamma \in \Gamma} \ p_{\HH^3}(\tilde{x},\gamma \cdot \tilde{y},t) \ .$$
\end{lemma}

The proof is rather straightforward. One shall first verify that the above right summation is finite which can be done using Formula \ref{formuleniyauchaleur}. Indeed, it shows that for any $t > 0$ the heat kernel $p_{\HH^3}$ decreases super-exponentially fast in $d(x,y)$. Next, one must verify that the summation is $\Gamma$-invariant with respect to both its spatial variables. It already holds (by construction) for the $y$-variable and it extends to the $x$-variable by symmetry. To conclude, now that the right member  of \eqref{formuleniyauchaleur} is well defined on $M_{\Gamma} \times M_{\Gamma} \times \RR_+^*$, it remains to verify that it solves the Cauchy problem \ref{eq Cauchy problem}, which is left to the reader. \\ 

Let us start the proof of Proposition \ref{propheatkernellip} by showing that the local version implies the compact version. \\

\textbf{Proof of (local version $\Rightarrow$ compact version).} Up to enlarging $K$, we will assume that $K$ is a closed ball containing both $x$ and $y$. We define 
	$$ U_n := \left\{ (w,z) \in K \times K \ , \ \limiinf{t \to \infty} \ p_{\Gamma}(w,z,t) \ t^{\alpha} > \frac{1}{n} \right\} \ , $$

and $$U := \underset{n \in \mathbb{N}} \cup U_n \ .$$

Since we suppose that $x,y \in K$ the set $U$ is non empty. Using the local version of Proposition \ref{propheatkernellip} and the symmetry of the heat kernel we get that all the sets $U_n$ are open, so is $U$. \\

We shall now see that $U$ is also closed. We fix $\epsilon = \frac{1}{2}$. Using again the local version of Proposition \ref{propheatkernellip} and the symmetry  we get that there is $\delta > 0$ such that for any points $w_1,z_1 \in \HH^3$ and any $C, \alpha > 0$ such that 
$$p_{\Gamma}(w_1,z_1,t) \ge C \ t^{-\alpha} \ ,$$
then there is a time $t_{\omega_1,z_1}$ such that for any $(w_2,z_2) \in B(w_1, \delta) \times B(z_1,\delta)$ and any $t > t_{\omega_1,z_1}$ we have 
	$$p_{\Gamma}(w_2,z_2,t) \ge \frac{C}{4} \ t^{-\alpha} \ ,$$

which shows that $U$ is closed. \\

We conclude by recalling that we supposed that $K$ is a ball and as such connected. This implies that $K \times K$ is connected too and therefore that $U =K \times K$. By definition of $U$ (as  an union of open sets) and because $K \times K$ is compact we get that there is a constant $C_K$ such that for any $w,z \in K$ 
	$$ \limiinf{t \to \infty} \ p_{\Gamma}(w,z,t) \ t^{\alpha} \ge C_K \ . $$ 

Using again the local version, one knows that for any pairs $(w_1,z_1) \in K \times K$ there is $\delta > 0$ and a time $t_{w_1,z_1} > 0$ such that for all $t \ge t_{w_1,z_1}$ and for all $(z_2,w_2) \in B(w_1,\delta) \times B(z_1,\delta)$ we have 
	$$ p_{\Gamma}(w_2,z_2,t) \ t^{\alpha} \ge \frac{C_K}{4} \ . $$

So that we covered $K \times K$ by open sets on which we have the above inequality uniformly in time. We conclude using again the compactness of $K \times K$ in order to get only finitely many such times and by taking the maximum of them. \hfill  $\blacksquare$ \\

\textbf{Proof of the local version of Proposition \ref{propheatkernellip}.} For what concerns this proof, we shall omit the tilde convention since we will only work on $\HH^3$. We start off splitting the summation appearing in Lemma \ref{lemheatkernelseries} according of how far an element $\gamma \in \Gamma$ moves $y$ away from $x$. Namely, for $\rho > 0$, we denote by 
	$$ 	B_{\Gamma}(x,y,\rho)  := \{ \ \gamma \in \Gamma \ , \ d(x, \gamma \cdot y ) \le \rho \ \} \ , $$
	
so that for any $z \in K$ we have 
	$$ p_{\Gamma}(x,z,t) = \somme{\gamma \in B_{\Gamma}(x,y,3t)} \ p_{\HH^3}(x,\gamma \cdot z,t) + \somme{\gamma \notin B_{\Gamma}(x,y,3t)} \ p_{\HH^3}(x,\gamma \cdot z,t)  \ .$$

We first show that the above right summation is negligible with respect to $p_{\Gamma}(x,y,t)$ when $t \to \infty$. More precisely, let us show that for every $\delta > 0$ and for all $z \in B(y, \delta)$ there is a constant $C(\delta,x,y) > 0$ such that for $t \ge 1$  
$$ \somme{\gamma \notin B_{\Gamma}(x,y,3t)} \ p_{\HH^3}(x,\gamma \cdot z,t)  \le C(\delta,x,y) \  e^{- \frac{t}{4}} \ . $$

\begin{figure}[!h]
	\begin{center}
		\def\svgwidth{0.6 \columnwidth}
			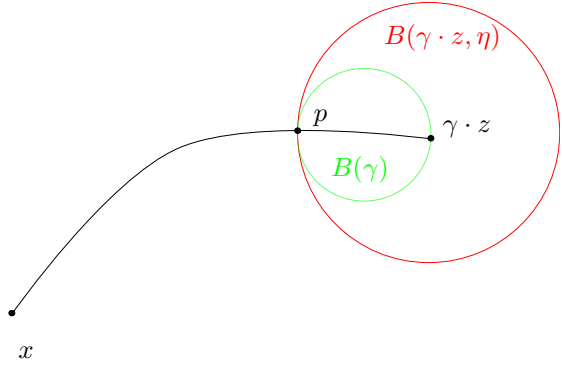
	\caption{The black curve represents the unique geodesic having endpoints $x$ and $\gamma \cdot z$. This geodesic intersects the ball $B(\gamma \cdot z, \eta)$ at $p$. The centre of the ball $B(\gamma)$ is chosen to be the midpoint of the geodesic of endpoints $\gamma \cdot z$ and $p$.} 
			\label{figure comptage plan}
	\end{center}
	\end{figure}

Since the ball $B(y, \delta)$ is compact, there is a constant $\eta$ such that for all $z \in B(y, \delta)$ the balls of $(B(\gamma \cdot z,\eta))_{\gamma \in \Gamma}$ are disjoints. Within any such a ball there is a ball of radius $\eta / 2$, denoted by $B(\gamma)$, such that for all $w \in B(\gamma)$ we have $d(x,\gamma \cdot z) \ge d(x,w) $ (see Figure \ref{figure comptage plan}). \\

Formula \ref{formuleniyauchaleur} shows that the function $ \rho \mapsto p_{\HH^3}(\rho,t) $ is decreasing. Therefore,
 $$p_{\HH^3}(x, \gamma \cdot z,t) \le \frac{1}{\vol(B(\gamma))}  \inte{B(\gamma)} p_\HH^3(x, w,  t) \ d \mu_h(w) \ .$$

Summing over $ \mathcal{B}_{\Gamma}(x,y,3t)$ we get
	$$\somme{\gamma \notin \mathcal{B}_{\Gamma}(x,y,3t) }  p_{\HH^3}(x, \gamma \cdot z,t) \prec  \somme{\gamma \notin \mathcal{B}_{\Gamma}(x,y,3t)}  \inte{ B(\gamma)} p_{\HH^3}(x, w,t) \ d \mu_h(w) \  , $$
where the symbol $\prec$ means that the left member is bounded by the right member up to a multiplicative constant which does not depend neither on $t$ nor on $z \in B(y, \delta)$
 .  \\

Note that all the balls $B(\gamma)$ are included in $ \{ w \in \HH^3 \ , \ d(x,w) > 3t - \delta - \eta \}$ and disjoints which yields
\begin{align*}
	 	\somme{\gamma \notin \mathcal{B}_{\Gamma}(x,y,3t)}  p_{\HH^3}(x, \gamma \cdot z,t) 
 		& \prec  \inte{ \{ d(x,w) \ge 3t - \delta - \eta  \} } p_{\HH^3}(x,w,t) \  d\mu_h(w)\\
 		& \prec \underset{{3t - \delta - \eta}}{\int^{+ \infty}} e^{2\rho}  p_{3}(\rho,t) \ d \rho \ ,
 \end{align*}
 
since the volume growth of balls in $\HH^3$ is bounded from above by $C \ e^{2\rho}$ for some $C > 0$. \\

Formula \eqref{formuleniyauchaleur} implies that for all $\rho \ge 1$ we have
	\begin{equation*}	
		 p_{3}(\rho,t) \prec \frac{\rho}{t^{\frac{3}{2}}}  e^{- t - \frac{\rho^2}{4t} - \rho} \ ,
	\end{equation*}	

and then 
		$$ e^{2  \rho}  p_{3}(\rho,t) \prec \frac{\rho}{t^{\frac{3}{2}}}  e^{ - \left( \frac{\rho - 2t}{2 \sqrt{t}} \right)^2} \ . $$

Therefore, for any $a > 1$ we have 
		$$ \inte{a}^{\infty} e^{2  \rho}  p_{3}(\rho,t) d \rho \prec \inte{a}^{\infty} \frac{\rho}{t^{\frac{3}{2}}} \ e^{ - \left( \frac{\rho - 2t}{2 \sqrt{t}} \right)^2} d \rho \ . $$

With the substitution $ u = \frac{\rho - 2t}{2 \sqrt{t}} $ we get 
	$$ \inte{a}^{\infty} \frac{\rho}{t^{\frac{3}{2}}} \ e^{ \left( \frac{\rho - 2t}{2 \sqrt{t}} \right)^2} d \rho = \inte{\frac{a - 2t}{2 \sqrt{t}}}^{\infty} \frac{ 2 \sqrt{t} u + 2t }{t^{\frac{3}{2}}} \ e^{ - u^2 } 2 \sqrt{t} \ d u \ .  $$

Setting $a = 3t - \delta - \eta$ we get for $t \ge 1$  and $u \ge \frac{a - 2t}{2 \sqrt{t}} \succ \sqrt{t} $ 
 $$ \frac{ 2 \sqrt{t} u + 2t }{t^{\frac{3}{2}}} \prec \frac{u}{\sqrt{t}} \ .$$ 
 Therefore, 
	\begin{align*}
		 \inte{[3t - \delta - \eta, \infty]} \frac{\rho}{t^{\frac{3}{2}}} \ e^{ \left( \frac{\rho - 2t}{2 \sqrt{t}} \right)^2} d \rho & \prec  \inte{[\frac{t - \delta - \eta}{2 \sqrt{t}}, \infty]} u \ e^{ - u^2 } d u \   \\
		 & \prec e^{  - \left( \frac{t - \delta - \eta}{2 \sqrt{t}} \right)^2 } \\
		 & \prec e^{  -  \frac{t }{4} } \ ,
	\end{align*}

and then
\begin{align*}
	 \somme{\gamma \notin \mathcal{B}_{\Gamma}(x,y,3t)}  p_{\HH^3}(x, \gamma \cdot z,t) \prec e^{ - \frac{t}{4} } \ .
\end{align*}

Because we supposed $p_{\Gamma}(x,y,t) \ge C \ t^{-\alpha}$ we have in particular
	$$ \somme{ \gamma \notin \mathcal{B}_{\Gamma}(x,y,3t) } p_{\HH^3}(x, \gamma \cdot z,t) \underset{t \to \infty}{=} o \Big( p_{\Gamma}(x,y,t)  \Big) ,  $$ 
uniformly on $z \in B(y,\delta)$. By uniform we mean that the underlying quantifiers of the above convergence do not depend on $z \in B(y, \delta)$. To sum up, we have proven so far that for any $z \in B(y, \delta)$ we have 
	$$ p_{\Gamma}(x,z,t) \underset{t \to \infty} =  \somme{ \gamma \in \mathcal{B}_{\Gamma}(x,y,3t) } p_{\HH^3}(x, \gamma \cdot z,t) + o \Big( p_{\Gamma}(x,y,t)  \Big) \ , $$ 

uniformly on $z \in B(y,\delta)$. In particular:
\begin{align*}
	 | p_{\Gamma}(x,z,t) - & p_{\Gamma}(x,y,t) | \underset{t \to \infty}{=}  \\ 
	 	 \Big| \somme{ \gamma \in \mathcal{B}_{\Gamma}(x,y,3t) } & p_{\HH^3}(x, \gamma \cdot z,t) - \somme{ \gamma \in \mathcal{B}_{\Gamma}(x,y,3t) } p_{\HH^3}(x, \gamma \cdot y,t) \Big| + o \Big( p_{\Gamma}(x,y,t)  \Big) .  
\end{align*}

It remains then to prove that the above absolute value is negligible with respect to $\ p_{\Gamma}(x,y,t)$ when $t \to \infty$ uniformly on $z \in B(y,\delta)$. \\

Namely, we shall see that for all $\epsilon > 0$ there is $\delta >  0$ such that for any $z \in B(y, \delta)$ we have 
	\begin{align*}
		 \Big| \somme{ \gamma \in \mathcal{B}_{\Gamma}(x,y,3t) } p_{\HH^3}(x, \gamma \cdot z,t) - \somme{ \gamma \in \mathcal{B}_{\Gamma}(x,y,3t) } p_{\HH^3}(x, \gamma \cdot y,t) \Big| & \le \epsilon   \somme{ \gamma \in \mathcal{B}_{\Gamma}(x,y,3t) } p_{\HH^3}(x, \gamma \cdot y,t) \ , 
	\end{align*}
which implies in particular
$$		\Big| \somme{ \gamma \in \mathcal{B}_{\Gamma}(x,y,3t) } p_{\HH^3}(x, \gamma \cdot z,t) - \somme{ \gamma \in \mathcal{B}_{\Gamma}(x,y,3t) } p_{\HH^3}(x, \gamma \cdot y,t) \Big| 	  \le \epsilon \ p_{\Gamma}(x,y,t) \ ,  $$
 concluding. \\

To do so, we will use again the explicit formula \ref{formuleniyauchaleur} to show that for any pairs $y,z$ with $d(y,z) \le \delta$ and for any $\gamma \in B_{\Gamma}(x,y,3t)$ we have for $t$ large enough that
 $$ | p_{\HH^3}(x, \gamma \cdot y,t) - p_{\HH^3}(x, \gamma \cdot z,t) | \le \epsilon \  p_{\HH^3}(x, \gamma \cdot y,t) \ , $$
which concludes the proof in summing over $\mathcal{B}_{\Gamma}(x,y,3t)$. \\

In fact, Formula \eqref{formuleniyauchaleur}  gives readily that there is two constants $C_1, C_2 > 0$ such that for all $x,y$ and for all $t > 1$ we have
	$$ \Big| \nabla_z \ln \Big( p_{\HH^3}(x,z,t) \Big) \Big| \le C_1 + C_2 \ \frac{\rho}{2t} \ ,  $$
with $\rho = d(x,z)$. So that, provided that $d(x,z) \le 4t$, we get the existence of $C_3 > 0$ such that
	$$ \Big| \nabla_z \ln \Big( p_{\HH^3}(x,z,t) \Big) \Big| \le 2 \ C_3 \ . $$
In particular, for any $z \in B(y, \delta)$ and for any $\gamma \in B_{\Gamma}(x,y,3t)$ we get 
	$$ \Big| \ln \left( \frac{p_{\HH^3}(x,\gamma \cdot z,t)}{p_{\HH^3}(x,\gamma \cdot y,t)} \right)  \Big| \le 2 \ C_3 \ \delta \ . $$

Therefore, for $\delta$ small enough
$$ \Big|  \frac{p_{\HH^3}(x,\gamma \cdot z,t)}{p_{\HH^3}(x,\gamma \cdot y,t)} - 1 \Big| \le \epsilon \ , $$

which is the desired conclusion. 	 \hfill $\blacksquare$

\section{Proof of Theorem \ref{maintheo}: reduction to the spectral proposition.}

\label{secproofmaintheo}

This Section is devoted to the reduction of the proof of our main theorem to Proposition 'spectral proposition'. We postpone the proof of the latter to the last section of this article. This section itself is split into two steps. We, most of the time, work in dimension 3 since we will use Proposition \ref{propheatkernellip} which requires it. \\

We shall first see how to deduce our main theorem from Proposition \ref{propmaintheo}. This part of the proof may be seen as how to recover the behaviour of the kernels from those of the operators. Afterwards, we will reduce the proof of Proposition 'operators comparison' to the proof of Proposition 'spectral proposition'. In order to soften the notations, we denote by $\left< \cdot, \cdot \right>$ the scalar product  $\left< \cdot, \cdot \right>_{L^2(M_{\Gamma})}$. Let us recall here the statement of Proposition \ref{propmaintheo} for the reader's convenience.

\begin{proposition}[operators comparison]
	Let $\Gamma$ be a Poincaré group such that there exist $x,y \in \mathbb{H}^3$ a and $\alpha > 0$ such that for $t$ large enough we have
		$$p_{\Gamma}(x,y,t) \ge t^{-\alpha} \ .$$ 
Then, for any non zero $f \in B^+_0(M_{\Gamma})$ 
		$$ \left< \mathcal{O}_{\Gamma}^{\rho}(f) , f \right> \equivaut{\rho \to \infty} \left< e^{- \frac{\rho}{2} \Delta}(f) , f \right> \ .$$
\end{proposition}

\textbf{Proof of (Proposition \ref{propmaintheo} $\Rightarrow$ Theorem \ref{maintheo}).} The following lemma, valid in any dimension $d$, is the geometric key of the proof. It is a variation around an argument of \cite{arteskinmac}. We set for any $x,y \in M_{\Gamma}$ and $\delta > 0$ 
$$ \text{\large{X}}_{\delta} := \frac{\mathds{1}_{B(x,\delta)}}{\mathcal{V}_d(\delta)} \hspace{1cm} \text{\large{Y}}_{\delta} := \frac{\mathds{1}_{B(y,\delta)}}{\mathcal{V}_d(\delta)} \ . $$
Recall that there is a constant $V_d$ such that $\mathcal{V}_d(\rho) \equivaut{\rho \to \infty} V_d \ e^{ 2 \rho}$. 

\begin{lemma}
	\label{lemmapproxorbital}
	Let $\Gamma$ be a Poincaré group, $x,y \in \mathbb{H}^d$. If there is a decreasing function $g : \RR_+ \to \RR_+$ such that for all $\epsilon > 0$ there is $\delta_{\epsilon} > 0$ such that for every $ 0 < \delta \le \delta_{\epsilon}$ and for $\rho$ large enough 
		$$ (1 - \epsilon) \ g(\rho) \le \left< \mathcal{O}_{\Gamma}^{\rho}(\text{\large{X}}_{\delta}) , \text{\large{Y}}_{\delta} \right> \le (1 + \epsilon) \ g(\rho) \ , $$
then 
	$$ N_{\Gamma}(x,y,\rho)  \equivaut{\rho \to \infty} g(\rho) \ \mathcal{V}_d(\rho) \equivaut{\rho \to \infty}  V_d \ g(\rho) \ e^{(d-1)\rho}. $$
\end{lemma}

\textbf{Proof of Lemma \ref{lemmapproxorbital}.} Given $\epsilon > 0$ we shall prove under the assumptions of Lemma \ref{lemmapproxorbital} that for $\rho$ large enough we have 
	$$  (1 -\epsilon)^2 \ g(\rho)  \le \frac{N_{\Gamma}(x,y,\rho)}{\mathcal{V}_d(\rho)}  \le (1 +\epsilon)^2 \ g(\rho) \ . $$
	
\begin{figure}[!h]
	\begin{center}
		\def\svgwidth{0.6 \columnwidth}
			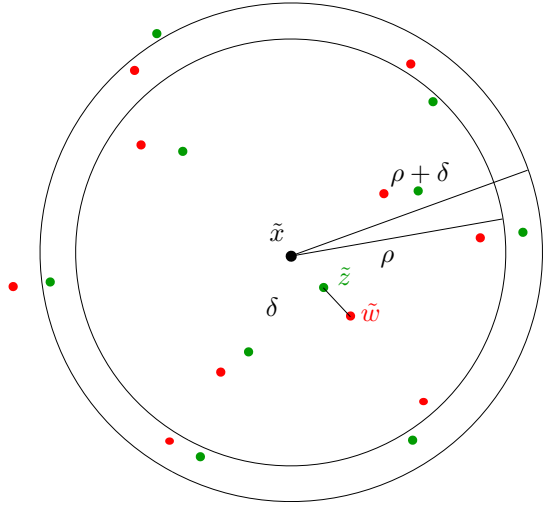
	\caption{the red dots correspond to the $\tilde{w}$-orbit  and the green ones to the $\tilde{z}$-orbit. Any red point can be paired with a green one by taking its closest neighbour (which are $\delta$-close from one another). Inequality \eqref{eqapproxorbital} corresponds to the fact that if a red point is in $B(\tilde{x},\rho)$ its green-mate must be in $B(\tilde{x}, \rho + \delta)$.}  
			\label{figure dessin}
	\end{center}
	\end{figure}

We start by working on $\HH^d$. The proof relies on the following key remark \cite[Section 2]{arteskinmac}, see Figure \ref{figure dessin}: for any $\delta >0$ and any $\tilde{x_0},\tilde{w},\tilde{z} \in \mathbb{H}^d$ such that $d(\tilde{w},\tilde{z}) \le  \delta$ we have 
	\begin{equation}
		\label{eqapproxorbital}
		 N_{\Gamma}(\tilde{x_0},\tilde{z}, \rho - \delta) \le N_{\Gamma}(\tilde{x_0},\tilde{w}, \rho) \le N_{\Gamma}(\tilde{x_0},\tilde{z}, \rho + \delta) \ . 
	\end{equation}
	
In particular, fixing $\tilde{x},\tilde{y} \in \mathbb{H}^d$, $ \tilde{\omega} \in B(\tilde{x},\epsilon)$ and $\tilde{z} \in B(\tilde{y},\epsilon)$ one has 
	$$ N_{\Gamma}(\tilde{w},\tilde{z}, \rho - 2 \delta ) \le N_{\Gamma}(\tilde{x},\tilde{y}, \rho) \le N_{\Gamma}(\tilde{w},\tilde{z}, \rho + 2 \delta) \ , $$
using twice \eqref{eqapproxorbital}. Averaging the above formula with respect to $\tilde{z}$ and $\tilde{w}$ on $B(\tilde{y},\delta) \times B(\tilde{x},\delta) $ and projecting to $M_{\Gamma}$ yields 
$$ \mathcal{V}_d(\rho - 2 \delta) \ \left< \mathcal{O}^{\Gamma}_{\rho - 2 \delta}(\text{\large{X}}_{\delta}) , \text{\large{Y}}_{\delta} \right> \le  N_{\Gamma}(x,y,\rho) \le \mathcal{V}_d(\rho + 2 \delta) \left< \mathcal{O}^{\Gamma}_{\rho + 2 \delta}(\text{\large{X}}_{\delta}) , \text{\large{Y}}_{\delta} \right> \ , $$

and then 
$$ \frac{\mathcal{V}_d(\rho - 2 \delta)}{\mathcal{V}_d(\rho)} \ \left< \mathcal{O}^{\Gamma}_{\rho - 2 \delta}(\text{\large{X}}_{\delta}) , \text{\large{Y}}_{\delta} \right> \le  \frac{N_{\Gamma}(x,y,\rho)}{\mathcal{V}_d(\rho)} \le \frac{\mathcal{V}_d(\rho + 2 \delta)}{\mathcal{V}_d(\rho)} \left< \mathcal{O}^{\Gamma}_{\rho + 2 \delta}(\text{\large{X}}_{\delta}) , \text{\large{Y}}_{\delta} \right> \ . $$

Using the assumption of Lemma \ref{lemmapproxorbital}, there is $\delta_{\epsilon} > 0$ such that for all $\delta_{\epsilon} >\delta > 0$ we have for $\rho$ large enough
	$$ \frac{\mathcal{V}_d(\rho - 2 \delta)}{\mathcal{V}_d(\rho)}  (1 -  \epsilon) g(\rho - 2 \delta) \le \frac{ N_{\Gamma}(x,y,\rho)}{\mathcal{V}_d(\rho)} \le \frac{\mathcal{V}_d(\rho + 2 \delta)}{\mathcal{V}_d(\rho)} (1 + \epsilon) g (\rho + 2 \delta)  \ . $$
Since we assume $g$ to be decreasing we also have 
$$ \frac{\mathcal{V}_d(\rho - 2 \delta)}{\mathcal{V}_d(\rho)}  (1 -  \epsilon) g(\rho) \le \frac{ N_{\Gamma}(x,y,\rho)}{\mathcal{V}_d(\rho)} \le \frac{\mathcal{V}_d(\rho + 2 \delta)}{\mathcal{V}_d(\rho)} (1 + \epsilon) g (\rho)  \ . $$

Moreover, since we know that there is a constant $V_d$ such that $\mathcal{V}_d(\rho) \equivaut{\rho \to \infty} V_d \ e^{(d-1) \rho}$, one gets that there is $\delta_0$ such that for $\delta \le \delta_0$ and $\rho$ large enough 
$$ 1 - \epsilon \le \frac{\mathcal{V}_d(\rho + 2 \delta)}{\mathcal{V}_d(\rho)} \le 1 + \epsilon \ . $$

Therefore, up to reducing $\delta_{\epsilon}$ such that $\delta_{\epsilon} \le \delta_0$, one gets that for all $\epsilon > 0$ and for $\rho$ large enough 
$$ (1 -  \epsilon)^2 g(\rho) \le \frac{ N_{\Gamma}(x,y,\rho)}{\mathcal{V}_d(\rho)} \le (1 + \epsilon)^2 g (\rho)  \ , $$

concluding. \hfill $\blacksquare$ \\

We now work with $d =3$. We want to use Lemma \ref{lemmapproxorbital} with $g(\rho) = p_{\Gamma}\left(x,x,\frac{\rho}{2}\right)$ together with Proposition \ref{propmaintheo}. To do so, we fix $\epsilon > 0$ and we use the local version of Proposition  \ref{propheatkernellip} which gives that there is $\delta_0 >0$ such that for $\delta \le \delta_0$ we have for $\rho$ large enough 
 $$ (1 - \epsilon) \ p_{\Gamma}\left(x,x,\frac{\rho}{2} \right) \le  \left< e^{- \frac{\rho}{2} \Delta}(\text{\large{X}}_{\delta} ) , \text{\large{X}}_{\delta}  \right> \le (1 + \epsilon) \ p_{\Gamma}\left(x,x,\frac{\rho}{2} \right) \ . $$

We now use Proposition \ref{propmaintheo} with $f =  \text{\large{X}}_{\delta} $. Combined with the above inequality it yields: for $\delta \le \delta_0$ we have for $\rho$ large enough
 $$ (1 - \epsilon)^2 \ p_{\Gamma}\left(x,x,\frac{\rho}{2} \right) \le \left< \mathcal{O}_{\Gamma}^{\rho}(\text{\large{X}}_{\delta} ) , \text{\large{X}}_{\delta} \right> \le (1 + \epsilon)^2 \ p_{\Gamma}\left(x,x,\frac{\rho}{2} \right) \ .$$

Lemma \ref{lemmeheatdecreases} asserts that $p_{\Gamma}\left(x,x,\frac{\rho}{2} \right)$ decreases as a function of $\rho$. Therefore, we have all the assumptions required to use Lemma \ref{lemmapproxorbital} which conclusion leads to the desired result. \hfill  $\blacksquare$ \\

It remains then to prove Proposition \ref{propmaintheo} which will occupy the rest of this article. We split the proof in two main steps. The first step we address is how to deduce Proposition \ref{propmaintheo} from the spectral proposition. We recall it here for the reader's convenience. 
\begin{proposition}[Spectral proposition]
	\label{propapproxvp}
For any $\beta >0$ there is a constant $C > 0$ such that for all $\epsilon > 0$ there is $\rho_0 > 0$ such that for all $\rho > \rho_0$ we have 
\begin{enumerate} 
	\item for all $ \lambda \le \frac{\beta \ln (\rho) }{\rho} $ 
			$$ | \nu_{\rho}(\lambda) - e^{-\lambda \frac{\rho}{d-1}}| \le \epsilon \ e^{ - \lambda \frac{\rho}{d-1}}  \ ;$$ 
	\item for all $ \lambda \ge \frac{\beta \ln (\rho) }{\rho}$ 
		$$  |\nu_{\rho}(\lambda)|  \le C \ \rho^{-\frac{\beta}{d-1}} \ . $$ 
\end{enumerate}
\end{proposition}

We will prove the spectral proposition in the last section of this article. Recall that we gave an overview of the proof of Proposition \ref{propmaintheo} in the introduction.\\

\textbf{Proof of (Proposition \ref{propapproxvp}
$\Rightarrow$ Proposition \ref{propmaintheo}).}
We fix a non zero $f \in B^+_0(M_{\Gamma})$ and $\epsilon>0$. \\

We want to show that there is $\rho_0$ such that for any $\rho \ge \rho_0$ we have 
\begin{equation}
	\label{equation11}
	 \Big| \left< \mathcal{O}_{\Gamma}^{\rho}(f) , f \right>  - \left< e^{- \frac{\rho}{2} \Delta}(f) , f \right>  \Big| \le \epsilon \left< e^{- \frac{\rho}{2} \Delta}(f) , f \right>  \ .
\end{equation}

Let us start by showing that there is $\rho_0 > 0$ such that for any $\rho > \rho_0$ the two following requirements are fulfilled. Recall that $\alpha \ge 0$ is the exponent given by the lower bound required on the heat kernel.
	\begin{itemize}
		\item the conclusion of Proposition \ref{propapproxvp} with $\beta = 2 (\alpha +1)$;
		\item and 
			$$ \frac{ C' \ ||f||^2_{L^2}}{\rho^{\alpha + 1}} \le \epsilon \left< e^{- \frac{\rho }{2} \Delta}(f) , f \right> \ ,$$
where $C' = \max(1,C)$ with $C = C(\beta) = C(\alpha, d)$, the same constant than the one given by the second item of Proposition \ref{propapproxvp}.
	\end{itemize}	 
There is nothing to prove (yet) for the first item since we postpone the proof of Proposition \ref{propapproxvp} to the next section. Let us clarify the second item. We will use the lower bound on the heat kernel together with the compact version of Proposition \ref{propheatkernellip}. In fact, we start off
\begin{align*}
 	\left< e^{- \frac{\rho }{2} \Delta}(f) , f \right> & = \inte{M_{\Gamma} \times M_{\Gamma}} p_{\Gamma} \left(x,y, \frac{\rho}{2} \right) \  f(x) \ f(y) \ d \mu_h(x) \ d \mu_h(y) \\
 	& = \inte{\supp{f} \times \supp{f}} p_{\Gamma} \left(x,y, \frac{\rho}{2} \right) \  f(x) \ f(y) \ d \mu_h(x) \ d \mu_h(y) \ . 
 \end{align*}
The compact version of Proposition \ref{propheatkernellip} with $K := \supp{f}$ gives that there is a constant $C(f) > 0$ such that for any $x,y \in \supp{f}$ and any $\rho$ large enough we have
$$ p_{\Gamma} \left(x,y, \frac{\rho}{2} \right) \ge C(f) \ \rho^{- \alpha} \ .$$
Therefore, because $f$ is non-negative, one has
\begin{align*}
	\left< e^{- \frac{\rho }{2} \Delta}(f) , f \right>  & \ge C(f)  \ \rho^{-\alpha} \inte{K \times K} f(x) \ f(y) \ d \mu_h(x) \ d \mu_h(y) \\
	& \ge C(f) \ \rho^{-\alpha} || f ||_{L^1}^2 \ .
\end{align*}
 In particular,
$$ \frac{ \left< e^{- \frac{\rho }{2} \Delta}(f) , f \right> }{\rho} \ge \frac{C_2(f)}{\rho^{\alpha +1}} \  $$
for some constant $C_2(f)$. This is equivalent to
$$ \frac{C' \ ||f||_{L^2}}{\rho^{\alpha +1}} \le   \frac{ C_3(f)}{\rho} \left< e^{- \frac{\rho }{2} \Delta}(f) , f \right> \ , $$
for some constant $C_3(f)$. It concludes letting $\rho_0$ such that 
$ \frac{ C_3(f)}{\rho_0}  \le \epsilon$. \\

We fix $\rho \ge \rho_0$ as above. We now show that for all domains $\Omega$ with smooth boundary which contains a $\rho$-neighbourhood of $\text{supp}(f)$ we have
\begin{equation}
	\label{equation12}
	 \Big| \left< \mathcal{O}_{\Gamma}^{\rho}(f) , f \right> - \left< e^{- \frac{\rho}{2} \Delta_{\Omega}}(f) , f \right> \Big| \le \epsilon \left< e^{- \frac{\rho}{2} \Delta_{\Omega}}(f) , f \right> \ . 
\end{equation}

This readily leads to \eqref{equation11} by taking an exhaustion $(\Omega_n)_{n \in \NN}$ of $M$ by such domains, letting $n \to \infty$ and recalling the conclusion of Theorem \ref{theododziuk} for the compact set $\{ \rho/2 \} \times \text{supp}(f) \times \text{supp}(f)$. \\

The support of $f$ being included in $\Omega$, one can expand it with respect to the $L^2(\Omega)$-eigenbasis given by Theorem \ref{theospectral}:
	$$ f = \somme{\lambda \in \mathrm{Sp}(\Omega)} \ \left<\Phi^{\Omega}_{\lambda},f \right> \Phi^{\Omega}_{\lambda} \ .$$ 

Since we assumed moreover that a $\rho$-neighbourhood of $\supp{f}$ is included in $\Omega$, one can also expand the function $\mathcal{O}^{\rho}_{\Gamma}(f)$:
\begin{equation*}
				 \mathcal{O}^{\rho}_{\Gamma}(f) = \somme{\lambda \in \mathrm{Sp}(\Omega)} \ \left<\Phi^{\Omega}_{\lambda}, \mathcal{O}^{\rho}_{\Gamma}(f) \right> \Phi^{\Omega}_{\lambda} \ .
\end{equation*}

Because orbital operators are symmetric we have
	$$ \left<\Phi^{\Omega}_{\lambda}, \mathcal{O}^{\rho}_{\Gamma}(f) \right> = \left< \mathcal{O}^{\rho}_{\Gamma}(\Phi^{\Omega}_{\lambda}), f \right> \ . $$

We now use Proposition \ref{propselbergtransform} to get for any $x \in \text{supp}(f)$
$$ \mathcal{O}_{\Gamma}^{\rho} \left( \Phi^{\Omega}_{\lambda} \right)(x) =  \nu_{\rho}(\lambda) \ \Phi^{\Omega}_{\lambda}(x) \ , $$

and then
	$$ \left<\Phi^{\Omega}_{\lambda}, \mathcal{O}^{\rho}_{\Gamma}(f) \right> = \nu_{\rho}(\lambda)  \left< \Phi^{\Omega}_{\lambda}, f \right> \ . $$

We use Plancherel's formula and expand the scalar product
$\left< \mathcal{O}_{\Gamma}^{\rho}(f),f \right>$ with respect to $(\Phi^{\Omega}_{\lambda})_{\lambda \in \text{Sp}(\Omega)}$ which yields
	$$ \left< \mathcal{O}_{\Gamma}^{\rho}(f),f \right> = \somme{\lambda \in \mathrm{Sp}(\Omega)} \ \nu_{\rho}(\lambda)  \left<\Phi^{\Omega}_{\lambda},f \right>^2 \ .$$ 

We now split the above right summation at low frequencies according to $\frac{\beta \ln (\rho)} {\rho}$: 
 \begin{equation}
 	\label{equation13}
 	  \somme{\lambda \in \mathrm{Sp}(\Omega)} \ \nu_{\rho}(\lambda)  \left<\Phi^{\Omega}_{\lambda},f \right>^2 = \somme{ \lambda \le \frac{\beta \ln (\rho)} {\rho}} \ \nu_{\rho}(\lambda)  \left<\Phi^{\Omega}_{\lambda},f \right>^2  + \somme{ \lambda > \frac{\beta \ln (\rho)} {\rho}} \ \nu_{\rho}(\lambda)  \left<\Phi^{\Omega}_{\lambda},f \right>^2 \ . 
\end{equation}

We will deal with these two summations independently. We start by showing that the right one is negligible. We set $\rho_0$ in a way to have access to Proposition \ref{propapproxvp} with $\beta = 2 (\alpha +1)$ for $\rho \ge \rho_0$. Therefore, using the second item of Proposition \ref{propapproxvp} gives a constant $C$ (which depends only on $\beta$) such that for $\lambda >  \frac{\beta \ln(\rho)}{\rho}$ 
\begin{align*}	
	| \nu_{\rho} (\lambda) | \le C \ \rho^{-\frac{\beta}{2}} \le  \frac{C}{\rho^{\alpha +1}} \ .
\end{align*}	 

Therefore
\begin{align*}
 \left| \somme{ \lambda > \frac{\beta \ln (\rho)} {\rho}} \ \nu_{\rho}(\lambda)  \left<\Phi^{\Omega}_{\lambda},f \right>^2 	\right|	& \le  \frac{ C}{\rho^{\alpha+1}} \  \somme{ \lambda > \frac{\beta \ln (\rho)} {\rho}} \  \left<\Phi^{\Omega}_{\lambda},f \right>^2  \\
 	& \le  \frac{C}{\rho^{\alpha+1}} \ ||f||_{L^2}^2 \ .
\end{align*}

Because of the way we set $\rho_0$ and since $\rho \ge \rho_0$ we get
\begin{equation}
	\label{eq20}
	\left| \somme{ \lambda > \frac{\beta \ln (\rho)} {\rho}} \ \nu_{\rho}(\lambda)  \left<\Phi^{\Omega}_{\lambda},f \right>^2 \right|	\le \epsilon \left< e^{- \frac{\rho}{2} \Delta_{\Omega}}(f) , f \right> \ .
\end{equation}

 We are now going to compare the first summation appearing in Equation \eqref{equation13} to $e^{- \frac{\rho }{2} \Delta_{\Omega}}(f)$ through the use of the first item of Proposition \ref{propapproxvp}. In fact, since $\left<\Phi^{\Omega}_{\lambda},f \right>^2 \ge 0$  we have
\begin{align*}
	 \Big| \somme{ \lambda \le \frac{\beta \ln (\rho)} {\rho}} \ \nu_{\rho}(\lambda)  \left<\Phi^{\Omega}_{\lambda},f \right>^2 - \somme{ \lambda \le \frac{\beta \ln (\rho)} {\rho}} & \ e^{-\lambda \frac{\rho}{2}}   \left<\Phi^{\Omega}_{\lambda},f \right>^2  \Big| \\
	 & \le \somme{ \lambda \le \frac{\beta \ln (\rho)} {\rho}} \ | e^{-\lambda \frac{\rho}{2}} - \nu_{\rho}(\lambda) |  \left<\Phi^{\Omega}_{\lambda},f \right>^2 
\end{align*}

Using the first item of Proposition \ref{propapproxvp} we get
\begin{align*}
	  \somme{ \lambda \le \frac{\beta \ln (\rho)} {\rho}} \ | e^{-\lambda \frac{\rho}{2}} - \nu_{\rho}(\lambda) |  \left<\Phi^{\Omega}_{\lambda},f \right>^2 & \le \epsilon \somme{ \lambda \le \frac{\beta \ln (\rho)} {\rho}} \ e^{-\lambda \frac{\rho}{2}}   \left<\Phi^{\Omega}_{\lambda},f \right>^2 \\ 
	  & \le  \epsilon \left< e^{- \frac{\rho}{2} \Delta_{\Omega}}(f) , f \right> \ ,
\end{align*}

since 
\begin{equation}
	\label{eqdemmaintheo3000}
	 \left< e^{- \frac{\rho}{2} \Delta_{\Omega}}(f) , f \right> = \somme{ \lambda \in \mathrm{Sp}(\Omega)} \ e^{-\lambda \frac{\rho}{2}}   \left<\Phi^{\Omega}_{\lambda},f \right>^2 \ .
\end{equation}

Combined with the upper bound \eqref{eq20} obtained for the second summation we get
\begin{equation} 
	\label{equation23}
	 \Big| \left< \mathcal{O}_{\Gamma}^{\rho}(f),f \right> - \somme{ \lambda \le \frac{\beta \ln (\rho)} {\rho}} \ e^{-\lambda \frac{\rho}{2}}   \left<\Phi^{\Omega}_{\lambda},f \right>^2  \Big|  \le 2 \ \epsilon  \left< e^{- \frac{\rho}{2} \Delta_{\Omega}}(f) , f \right> \ .  
\end{equation}
	
It remains then to relate the above summation with the heat operator on $M_\Gamma$. This will be handle using again the choices made on $\beta$ and $\rho_0$. Indeed, since $\beta = 2(\alpha +1)$ we have for all $\lambda > \frac{\beta \ln (\rho)}{\rho}$
	$$ e^{-\lambda \frac{\rho}{2}}  \le \frac{1}{\rho^{\alpha +1}} \ . $$
Therefore, summing over $\lambda > \frac{\beta \ln (\rho)} {\rho}$ gives
$$ \somme{ \lambda > \frac{\beta \ln (\rho)} {\rho}} \ e^{-\lambda \frac{\rho}{2}} \left<\Phi^{\Omega}_{\lambda},f \right>^2 \le  \frac{||f||_{L^2}^2}{\rho^{\alpha+1}} \le \epsilon \left< e^{- \frac{\rho}{2} \Delta_{\Omega}}(f) , f \right> \ , $$
because of how we set $\rho_0$. Therefore, using \eqref{eqdemmaintheo3000} and the above inequality, we get
 $$ \Big| \left< e^{- \frac{\rho}{2} \Delta_{\Omega}}(f) , f \right> - \somme{ \lambda \le \frac{\beta \ln (\rho)} {\rho}} \ e^{-\lambda \frac{\rho}{2}}   \left<\Phi^{\Omega}_{\lambda},f \right>^2  \Big|  \le  \epsilon \ \left< e^{- \frac{\rho}{2} \Delta_{\Omega}}(f) , f \right> \ , $$

which, combined with \eqref{equation23} yields 
	$$ \Big| \left< \mathcal{O}_{\Gamma}^{\rho}(f),f \right> - \left< e^{- \frac{\rho}{2} \Delta_{\Omega}}(f) , f \right> \Big| \le 3 \ \epsilon \ \left< e^{- \frac{\rho}{2} \Delta_{\Omega}}(f) , f \right> \ ,  $$
concluding. \hfill $\blacksquare$

\section{Proof of the spectral proposition \ref{propapproxvp}}
\label{secetudevpselberg}

This section is devoted to the proof of Proposition \ref{propapproxvp}. It is quite technical and to ease the exposition we split it in three parts. We shall use a lot that $s$ is a solution of the equation $s(d-1-s) = \lambda$ with real part less that $(d-1)/2$. Since we shall also think of $s$ as a function of $\lambda$, we will often denote $s$ by $s(\lambda)$ when we think it clarifies the exposition (especially in Step 2). Note that 
	$$ \fonctionbis{[0, (d-1)^2/4 ]}{[0, (d -1)/2]}{\lambda}{s(\lambda)}$$
is an increasing diffeomorphism. 

\subsection{Proof of the spectral proposition: step 1.} In order to state our first step, we set for $\rho > |r|$.
	$$  \theta(\rho,r) := \left(\frac{2(\cosh(\rho) - \cosh(r))}{e^{\rho}} \right)^{\frac{d-1}{2}} \ . $$

Note that	for all $\rho > |r|$ 
	$$  0 \le \theta(\rho,r) \le 2^{\frac{d-1}{2}} \ .$$

\begin{lemma}
	\label{lemmaselbergstoformula}
For any $\lambda, \rho > 0$ we have : 
\begin{equation}
		\label{eqlemmeapproxvp}	 \mathcal{V}_d(\rho) \ \nu_{\rho}(\lambda) =  c_{d-1} e^{-\rho s} e^{(d-1)\rho} \inte{0}^{2\rho} \ e^{u \left(s - \frac{d-1}{2} \right)}  \theta(\rho, u-\rho)  \ du \ , 
\end{equation}

where $s = s(\lambda)$ is the unique complex number satisfying $s(d-1-s) = \lambda$ with $\Re(s) \le \frac{d-1}{2}$ and $c_{d-1}$ is volume of an euclidean ball in $\RR^{d-1}$ of radius $1$. 
\end{lemma}

\begin{figure}[!h]
	\begin{center}
		\def\svgwidth{0.6 \columnwidth}
			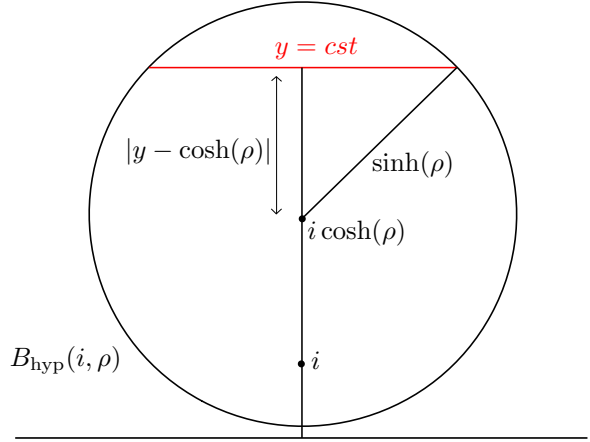
	\caption{The red part corresponds to an euclidean ball of dimension $d-1$.  The euclidean radius of such a ball is obtained by using the Pythagorean theorem.}  
			\label{figure dessin1}
	\end{center}
	\end{figure}

\textbf{Proof.} Let us denote by $I$ the left member of Equation \eqref{eqlemmeapproxvp} which, using Proposition \ref{propselbergtransform}, satisfies 
	\begin{equation}
		\label{eqfig1}
			 I = \inte{B(o,\rho)} \ y^s \  \frac{dx_1...dx_{d-1}dy}{y^d} \ , 
	\end{equation}
	
using the upper half model of the hyperbolic space $\HH^d$. In this model, an hyperbolic ball of radius $\rho$ centred at $(0,...,0,1)$ corresponds to an euclidean ball centred at $(0,....,0,\cosh(\rho))$ of radius $\sinh(\rho)$. Using Fubini's theorem we integrate first with respect to the $d-1$ first coordinates and then with respect to the last one (see Figure \ref{figure dessin1}) to get
$$ I = \inte{e^{-\rho}}^{e^{\rho}} \ y^s \ \vol_{\mathrm{euc}} \left( B \left( \sqrt{ \sinh(\rho)^2 - (y - \cosh(\rho))^2} \right) \right) \ \frac{dy}{y^d} \ , $$

where we denoted by $\vol_{\text{euc}} \left( B(r) \right)$ the euclidean volume of any ball of $\RR^{d-1}$ of radius $r$. \\

Recall that there is a constant $c_{d-1}$ such that 
	$$ \vol_{\text{euc}} \left( B(r) \right) = c_{d-1} r^{d-1} \ . $$
	
Which, combined with the identity $\sinh(\rho)^2 - \cosh(\rho)^2 = -1$ gives
	$$ I = \inte{e^{-\rho}}^{e^{\rho}} \ y^s \ c_{d-1} \left(  2 \cosh(\rho)y - y^2 - 1 \right)^{\frac{d-1}{2}}  \ \frac{dy}{y^d} \ . $$

With the substitution $y = e^r$ we obtain
	$$ I = c_{d-1} \inte{-\rho}^{\rho} \ e^{rs} \ \left(  2 \cosh(\rho)e^r - e^{2r} - 1 \right)^{\frac{d-1}{2}}  \ \frac{e^{r}dr}{e^{rd}} \ . $$

We shall prefer the following writing 
 $$ \frac{e^{r}}{e^{rd}} = e^{-(d-1)r} = e^{-\frac{(d-1)r}{2}} \cdot e^{-\frac{(d-1)r}{2}} $$
 
in order to get 
\begin{align*}
 I & = c_{d-1} \inte{-\rho}^{\rho} \ e^{rs} e^{-\frac{(d-1)r}{2}} \ \left( 2 \cosh(\rho) - ( e^{r} + e^{-r}) \right)^{\frac{d-1}{2}} \ dr \\ 
 	& = c_{d-1} \inte{-\rho}^{\rho} \ e^{r \left(s- \frac{d-1}{2} \right)} \ \left( 2 (\cosh(\rho) -  \cosh{r}) \right)^{\frac{d-1}{2}} \ dr \ . \\
 	& = c_{d-1} \inte{-\rho}^{\rho} \ e^{r \left(s- \frac{d-1}{2} \right)} \ \theta(\rho,r) \ e^{\rho \left(\frac{d-1}{2} \right) } \ dr \\ 
	& = c_{d-1} \inte{-\rho}^{\rho} \ e^{(-\rho) \left(s- \frac{d-1}{2} \right)} e^{(r + \rho) \left(s- \frac{d-1}{2} \right)} \ \theta(\rho,r) \ e^{\rho \left(\frac{d-1}{2} \right) } \ dr \\ 
 	& = c_{d-1} e^{(d-1)\rho} e^{ -s \rho} \inte{-\rho}^{\rho} \ e^{(r + \rho) \left(s- \frac{d-1}{2} \right)} \ \theta(\rho,r)  \ dr \ ,
\end{align*}	

which concludes the proof, substituing $r + \rho$ by $ u$ in the above integral. \hfill $\blacksquare$ 

\subsection{Proof of the spectral proposition: step 2.} 
We will now use the explicit formula given by Lemma \ref{lemmaselbergstoformula} to make our second step toward the spectral proposition \ref{propapproxvp}. Recall that $s(\lambda)$ is such as $s(\lambda)(d-1-s(\lambda)) = \lambda$.

\begin{lemma} 
\label{lemmeapproxvp}
Given any $ 0 < \lambda_0 < \left(\frac{d-1}{2} \right)^2 $ there is a constant $C = C(\lambda_0) > 0$ such that for $\rho \ge 1$ we have: 
	\begin{enumerate}
		\item If $\lambda \le \lambda_0$ then 
			$$ | \nu_{\rho}(\lambda) e^{s(\lambda) \rho} - 1 | \le C \ s(\lambda) \ ;$$ 		
		\item If $\lambda > \lambda_0$ then
	$$ |\nu_{\rho}(\lambda)| \le C \ \rho \ e^{- \rho s(\lambda_0)} \ . $$
	\end{enumerate}
\end{lemma}

\textbf{Proof.} First, we set 
$$ \mathcal{I}(s,\rho) := \inte{0}^{2\rho} \ e^{u \left(s - \frac{d-1}{2} \right)}  \theta(\rho, u-\rho)   \ du  $$

in such a way that \eqref{eqlemmeapproxvp} gives
$$ e^{s(\lambda) \rho} \nu_{\rho}(\lambda) = \frac{ c_{d-1} e^{(d-1)\rho}} {\mathcal{V}_d(\rho)} \ \mathcal{I}(s(\lambda),\rho) \ . $$

We start by working the second item of Lemma \ref{lemmeapproxvp} which is the easier to obtain. Let $ 0 < \lambda_0 < \left( \frac{d-1}{2} \right)^2 $  and $\lambda > \lambda_0$. We compute
\begin{align*}
	 |\nu_{\rho}(\lambda)| & \le \frac{c_{d-1} e^{(d-1)\rho}}{\mathcal{V}_d(\rho)} | e^{ - s(\lambda) \rho } \ \mathcal{I}(s(\lambda),\rho) | \\
	 	& \le \frac{c_{d-1} e^{(d-1)\rho}}{\mathcal{V}_d(\rho)}  e^{ - \Re(s(\lambda)) \rho } \ | \mathcal{I}(s(\lambda),\rho) | \ .
\end{align*}
Recall that $s(\lambda)$ is a solution of $s(d-1 -s) = \lambda$ of real part less that $1/2$. \\

Since $\frac{c_{d-1} e^{(d-1)\rho}}{\mathcal{V}_d(\rho)}$ is bounded from above by a constant $C_2$ we have for $\rho \ge 1$ 
	 	$$	| \nu_{\rho}(\lambda)|  \le C_2 \ e^{ - \Re(s(\lambda)) \rho } \ | \mathcal{I}(s(\lambda),\rho) | \ . $$

The function $ \lambda \mapsto \Re(s(\lambda))$ being increasing on the interval $[0,\left(\frac{d-1}{2}\right)^2]$ and constant equal to $\frac{d-1}{2}$ on $[\left(\frac{d-1}{2}\right)^2,\infty)$, we have 
	 	$$	| \nu_{\rho}(\lambda)|  \le C_2 \ e^{ - \Re(s(\lambda_0)) \rho } \ | \mathcal{I}(s(\lambda),\rho) | \ . $$

Moreover, since $\lambda_0 < \left(\frac{d-1}{2}\right)^2$ we actually have $\Re(s(\lambda_0)) = s(\lambda_0)$ which gives 
	 	$$	 |\nu_{\rho}(\lambda)|  \le C_2 \ e^{ - s(\lambda_0) \rho } \ | \mathcal{I}(s(\lambda),\rho) | \ . $$
We conclude for what concerns the second item by recalling that the function $\theta$ is bounded from above by $2^{ \frac{d-1}{2}}$. Therefore
\begin{align*}
	 | \mathcal{I}(s(\lambda),\rho) | & \le 2^{ \frac{d-1}{2}} \ \inte{0}^{2\rho} \  | e^{u \left(s - \frac{d-1}{2} \right)} | \ du \\
	 	 & \le 2^{ \frac{d-1}{2}} \ \inte{0}^{2\rho} \   e^{u \ \left( \Re(s(\lambda)) - \frac{d-1}{2} \right)} \ du \\
	 	 & \le 2^{ \frac{d+1}{2}} \rho  \ ,
\end{align*}

since $\Re(s(\lambda)) \le \frac{d-1}{2}$. Which yields
	$$ |\nu_{\rho}(\lambda)| \le  2^{ \frac{d+1}{2}} \rho \ C_2  \ e^{ - s(\lambda_0) \rho} \ , $$

the desired upper bound. \\

In order to prove that the first item of Lemma \ref{lemmeapproxvp} holds note that for any $\rho > 0$ we have $ \nu_{\rho}(0) = 1$, indeed:
	$$ \nu_{\rho}(0) = \finte{B(o, \rho)}{ y^{s(0)}}{  \frac{dx_1...dx_{d-1} dy}{y^d}}  = \finte{B(o, \rho)}{1}{  \frac{dx_1...dx_{d-1} dy}{y^d}}  = 1 \ ,$$
since $s(0) = 0$. Therefore $\nu_{\rho}(0) e^{s(0) \rho} = 1$ and for any $\rho \ge 1$ and any $\lambda \le \lambda_0$ we have
\begin{align*}
	\Big| \nu_{\rho}(\lambda) e^{s(\lambda) \rho} - 1 \Big| & = \frac{ c_{d-1} e^{(d-1)\rho}} {\mathcal{V}_d(\rho)} \ \Big| \mathcal{I}(s(\lambda),\rho) - \mathcal{I}(0,\rho) \Big| \\
	& \le C_2 \Big| \mathcal{I}(s(\lambda),\rho) - \mathcal{I}(0,\rho) \Big|  \\
	& \le C_2 \inte{0}^{s(\lambda)} \  \Big| \partial_1 \mathcal{I}(w,\rho) \Big| \ dw \ ,  
\end{align*}

since the function $\mathcal{I}$ is smooth and because $s(\lambda) \in \RR_+$ (since $\lambda_0 < (d-1)/2$). Differentiating under the integral and using again that the function $\theta$ is bounded we get that for any $ 0 \le w \le s(\lambda_0)$ 
\begin{align*}
 	\Big| \partial_1 \mathcal{I}(w,\rho) \Big| & \le 2^{ \frac{d-1}{2}} \inte{\RR_+} u \ e^{ u \ \left(w -\frac{d-1}{2} \right) } \ du \\
		 & \le 2^{ \frac{d-1}{2}} \inte{\RR_+} u \ e^{ u \ \left(s(\lambda_0) -\frac{d-1}{2} \right) } \ du \ ,
\end{align*}

using again that the function $\lambda \mapsto s(\lambda)$ is increasing on $[0, \lambda_0]$. Therefore, since we chose $\lambda_0$ such that $s(\lambda_0) < \frac{d-1}{2}$, we get a constant $C_3 = C_3(\lambda_0)$ such that for all $\rho > 1$ and all $\omega \le s(\lambda_0)$ we have
 	$$ \Big| \partial_1 \mathcal{I}(w,\rho) \Big|  \le C_3 \ , $$

and then
$$	\Big| \nu_{\rho}(\lambda) e^{s(\lambda) \rho} - 1 \Big|  \le C_2 \ C_3 \ s(\lambda) \ ,$$

which is the desired conclusion. \hfill $\blacksquare$ 

\subsection{Proof of the spectral proposition: step 3.}

We conclude by showing how Lemma \ref{lemmeapproxvp} implies the spectral proposition \ref{propapproxvp}. We shall start with the first point of Proposition \ref{propapproxvp}. Let $\rho \ge 1$ and set
 $$\varphi(\lambda,\rho) := | \nu_{\rho}(\lambda) - e^{-\lambda \frac{\rho}{d-1}}| \ .$$ 

Multiplying both sides by $e^{s \rho}$ gives 
	$$  e^{s \rho} \varphi(\lambda,\rho) = | \nu_{\rho}(\lambda)e^{s\rho} - e^{\left( s - \frac{\lambda}{d-1} \right) \rho}| \ . $$ 

Using the triangular identity and Lemma \ref{lemmeapproxvp} we get 
\begin{equation}
	\label{equation13bis}
	 e^{s \rho} \varphi(\lambda,\rho) \le C s + | 1 - e^{\left( s - \frac{\lambda}{d-1} \right) \rho} | 
\end{equation}

Recall that $s$ is defined as solving $s( d-1 -s) = \lambda$ with $\Re(s) \le \frac{d-1}{2}$. In particular we have 
\begin{equation}
	\label{equation14}
	\lambda - (d-1)s = -s^2 \ , 
\end{equation} 

from which we deduce $(d-1)s \equivaut{\lambda \to 0} \lambda$. Therefore, for $\rho$ large enough we have for all $\lambda \le \frac{\beta \ln(\rho) }{\rho}$
			$$ (d-1) \ s \le 2 \lambda \le 2 \  \frac{\beta \ln (\rho)}{\rho} \ ,$$
which gives
\begin{equation}
	\label{upperbounds}
		s \le \frac{2}{d-1} \frac{\beta \ln(\rho)}{\rho} \ . 
\end{equation}
			
Looking backward to Equation \eqref{equation14}, we get a constant $C_2 > 0$ such that
\begin{equation*}
			0 \le \left( s - \frac{\lambda}{d-1} \right) \rho = \frac{s^2}{d-1} \le C_2 \frac{ \ln^2(\rho)}{\rho} \ .  \\
\end{equation*} 

The above identity combined with \eqref{equation13bis} and the upper bound \eqref{upperbounds} gives for $\rho$ large enough 
	$$ e^{s \rho} \varphi(\lambda,\rho) \le C_2 \cdot \frac{2 \beta}{d-1} \frac{\ln(\rho)}{\rho} + e^{ C_2 \frac{ \ln^2(\rho)}{\rho}} -1 \ . $$
	
Because of $\frac{\ln^2(\rho)}{\rho} \tends{\rho \to \infty} 0$, one has for $\rho$ large enough that 
	$$ e^{ C_2 \frac{ \ln^2(\rho)}{\rho}} -1 \le 2 \ C_3 \ \frac{\ln^2(\rho)}{\rho} \ .$$ 

Therefore, there are two constants $C_3, C_4$ such that for $\rho$ large enough 
\begin{align*}
	 e^{s \rho} \varphi(\lambda,\rho) & \le C_3 \cdot \left( \frac{\ln(\rho)}{\rho} +  \frac{ \ln^2(\rho)}{\rho} \right) \ , \\
	 	& \le C_4 \frac{ \ln^2(\rho)}{\rho} 
\end{align*}	

and 
 $$  \varphi(\lambda,\rho)  \le C_4 \frac{ \ln^2(\rho)}{\rho} \ e^{- s \rho} \ . $$

Recall also that $ s \ge \frac{\lambda}{d-1}$ which yields 
	$$  \varphi(\lambda,\rho) \le C_4 \frac{ \ln^2(\rho)}{\rho} e^{ - \lambda \frac{\rho}{d-1}} \ \ . $$
This concludes by setting $\rho$ large enough in order for 
	$$ C_4 \frac{ \ln^2(\rho)}{\rho}  \le \epsilon $$ 
to hold. \\

We now deal with the second item of Proposition \ref{propapproxvp}. Let any  $\lambda_0$ such that $0 < \lambda_0 < (d-1)/2$ and $\rho_0$ as in Lemma \ref{lemmeapproxvp} for this given $\lambda_0$. We split the proof depending on whether $\lambda \le \lambda_0$ or not. Lemma \ref{lemmeapproxvp} asserts that for all $\rho \ge \rho_0$ and for any $ \frac{\beta \ln \rho}{\rho} \le \lambda \le \lambda_0$ we have 
\begin{align*}
	\nu_{\rho}(\lambda)  e^{s \rho} & \le 1 + C \ s \\
		& \le C_6 \ , 
\end{align*}

with $C_6 : = \underset{\lambda \le \lambda_0}{\sup} \{ 1 + C s(\lambda) \}$. Hence, 
	$$ \nu_{\rho}(\lambda) \le C_6 \ e^{-s \rho}   \ .  $$

Thus, using again $ \frac{\lambda}{d-1} \le s$, we get
	$$ \nu_{\rho}(\lambda) \le C_6 \ e^{- \frac{\lambda}{d-1} \rho}   \ .  $$

But, since we assumed $\lambda \ge \frac{\beta  \ln(\rho) }{\rho}$, we also have 
	$$ \nu_{\rho}(\lambda) \le C_6 \ e^{- \frac{\beta}{d-1} \ln(\rho)}   \ ,  $$

which gives the expected result for $\lambda \le \lambda_0$.  If $\lambda > \lambda_0$ we use the second part of Lemma \ref{lemmeapproxvp} which implies that the function $|\nu_{\rho}(\lambda)|$ decreases exponentially fast to $0$. In particular, faster than $\rho^{-\frac{\beta}{d-1}}$. \hfill $\blacksquare$ 

\appendix

\section{Eigenvectors of Radial operators}
\label{secselbergspretrace}
This appendix aims at proving Proposition \ref{propselbergtransform} which is the first key observation toward Selberg's trace formula. It takes its root back to Delsarte's note \cite{artdelsarte}. It is classical, but we decided to prove it in this appendix since there is no ready-to-use statement in the literature for open domains (to the author's best knowledge). Note that the following discussion is valid for any symmetric space of rank $1$ and any radial operators, meaning any operator whose kernel $\kappa(x,y)$ is given by a function of $d(x,y)$, see \cite[Chapter 11]{livrechavel} for more details. Let us recall Proposition \ref{propselbergtransform} for the reader's convenience. We use the upper half plane model of the hyperbolic space $\mathbb{H}^d$ which is given by $\RR^{d-1} \times \RR_+$ endowed with the metric 
	$$ \frac{dx^2_1 + ... + dx^2_{d-1} + dy^2}{y^{2}} \ . $$

\begin{proposition}[Delsarte's formula]
Let $\Gamma$ be a Poincaré group and $\Omega$ any open subset of $M_{\Gamma}$ and $\Phi$ a smooth function of $\Omega$ verifying $\Delta \Phi = \lambda \Phi$. Then for any $x \in \Omega$ and any $\rho > 0$ such that $B(x,\rho) \subset \Omega$ we have
\begin{equation*}
	 \mathcal{O}_{\Gamma}^{\rho}(\Phi)(x) = \nu_{\rho}(\lambda) \ \Phi(x) \ ,
\end{equation*}

with 
\begin{equation*}
\nu_{\rho}(\lambda) = \mathcal{O}^{\rho}(y^s) = \finte{B(o,\rho)} { y^s}{  \frac{dx_1...dx_{d-1} dy}{y^d}}  
\end{equation*}
where $o := (0,..., 0, 1) \in \RR^{d-1} \times \RR_+^*$ and $s$ satisfying $s(d-1-s) = \lambda$.
\end{proposition}

Let $\Omega, \Phi, x$ and $\rho$ as in the above proposition. By definition
 $$ \mathcal{O}^{\rho}_{\Gamma}(\Phi)(x) := \frac{1}{\mathcal{V}_d(\rho)} \inte{B(\tilde{x},\rho)}  \tilde{\Phi}(y) \ d \mu_h(y) \ ,  $$

where $\tilde{\Phi}$ is the lift of the function $\Phi$ on $\pi^{-1}(\Omega) \subset \HH^d$. Because $\pi$ is a local isometry, the function $\tilde{\Phi}$ solves $\Delta = \lambda \Id$ on $\pi^{-1}(\Omega)$ which is also an open domain of $\HH^d$. Therefore, Proposition \ref{propselbergtransform} follows from the following proposition (which corresponds to the special case $\Gamma = \{\Id\}$).

\begin{proposition}
	\label{propselbergprepretrace}
Let $\Phi$ be a smooth function of an open set $\Omega \subset \HH^d$ verifying $\Delta \Phi = \lambda \Phi$. Then for any $x \in \Omega$ and any $\rho$ such that $B(x,\rho) \subset \Omega$ we have
\begin{equation*}
	 \mathcal{O}^{\rho}(\Phi)(x) = \nu_{\rho}(\lambda) \ \Phi(x) \ ,
\end{equation*}
with 
\begin{equation*}
	\nu_{\rho}(\lambda) = \mathcal{O}^{\rho}(y^s) = \finte{B(o,\rho)} { y^s}{  \frac{dx_1...dx_{d-1} dy}{y^d}}  
\end{equation*}
where $o := (0,..., 0, 1) \in \RR^{d-1} \times \RR_+^*$ and $s$ satisfying $s(d-1-s) = \lambda$.
\end{proposition}

The end of this appendix is devoted to the proof of the above statement. The key property used in order to establish Proposition \ref{propselbergprepretrace} is that radial functions of symmetric spaces of rank 1 satisfying the partial differential equation $\Delta = \lambda \Id$ are solutions of an ordinary differential equation (abbreviated ODE). \\

\textbf{Proof.} The following lemma is the intermediate key. We denote by $S(x,\rho)$ the sphere centred at $x$ of radius $\rho$ and by $d \mu_{S(x,\rho)}$ the measure on $S(x,\rho)$ induced by the hyperbolic metric.

\begin{lemma}
	\label{lemmaappendix2}
			There is an 'universal'  function $\mathfrak{S} : \RR_+ \times \RR_+^* \to \RR$ with the property that that for any $\rho >0$, $\lambda \ge 0$, any open set $\Omega \subset \HH^d$ and any smooth function $\Phi : \Omega \to \CC$ satisfying $\Delta \Phi = \lambda \Phi$ we have for all $x \in \Omega$ with $B(x,\rho) \subset \Omega$ 
	$$ \finte{S(x, \rho)}{\Phi}{ d \mu_{S(x,\rho)}} = \mathfrak{S}(\lambda,\rho) \ \Phi(x) \ .$$
\end{lemma}

Before proving the above lemma, let us show how it implies Proposition \ref{propselbergprepretrace}. \\

\textbf{Proof of (Lemma }\ref{lemmaappendix2} $\Rightarrow$ \textbf{Proposition \ref{propselbergprepretrace}).} We start by showing that the function $\nu_{\rho}$ are well defined, which corresponds to the first part of Proposition \ref{propselbergprepretrace}. Let $\Omega, \Phi, x$ and $\rho$ as above. Since the volume of $S(x, \rho)$ is given by $\rho \mapsto \mathcal{V}_d'(\rho)$, one has
\begin{align*}
	 \finte{B(x,\rho)} \Phi \ d \mu_h  = \frac{1}{\mathcal{V}_d(\rho)} \inte{0}^{\rho}
\mathcal{V}_d'(t) \left(\finte{S(x, t)}{\Phi}{ d \mu_{S(x,t)}} \right) dt 	 \ .
\end{align*}
Using Lemma \ref{lemmaappendix2} we get
\begin{equation}
\label{eqappend1}
   \finte{B(x,\rho)}{\Phi}{d \mu_h} = \left( \frac{1}{\mathcal{V}_d(\rho)} \inte{0}^{\rho}
\mathcal{V}_d'(t) \ \mathfrak{S}(\lambda,t) \ dt \right)  \Phi(x) \ ,
\end{equation}
which concludes setting  $$\nu_{\rho}(\lambda) := \frac{1}{\mathcal{V}_d(\rho)} \inte{0}^{\rho}
\mathcal{V}_d'(t) \ \mathfrak{S}(\lambda,t) \ dt \ . $$

We now prove the second part of Proposition  \ref{propselbergprepretrace} which gives an explicit formula for the Selberg's transform. It is an easy consequence of the existence of $\nu_\rho$. \\

Given $\lambda > 0$ let $s$ be a solution of 
	$$ \lambda = s (d-1-s) \ . $$ 

Using the following general expression of the Laplace operator in coordinates (we use Einstein's summation convention),
$$\Delta = \frac{1}{\sqrt{\det(g_{ij})}} \partial_i \left( \sqrt{\det(g_{ij})} g^{ij} \partial_j \right),$$

 one easily computes that for any $s \in \CC$ the functions
	$$ \fonction{\Psi_{s}}{\RR^{d-1} \times \RR_+}{\CC}{(x_1, ... , x_{d-1}, y)}{y^s} $$
	
satisfies $\Delta \Psi_{s} = \lambda \  \Psi_{s}$ with $\lambda = s(d-1-s)$ and values $1$ at the point $o := (0,...,0, 1)$. \\ 

Because of Equation \ref{eqappend1} and because of our definition of $\nu_{\rho}$ we have
	$$ \nu_{\rho}(\lambda) = \nu_{\rho}(\lambda) \ \Psi_{s}(o) = \underset{{B(o,\rho)}}{\fint}  \Psi_{s} \ d \mu_h = \finte{B(o,\rho)} { y^s}{  \frac{dx_1...dx_{d-1} dy}{y^d}} \ ,  $$
concluding.  \hfill $\blacksquare$ \\

\textbf{Proof of Lemma \ref{lemmaappendix2}.} The first step toward the proof is to reduce the study to radial functions. Given a point $x \in \Omega$, we say that a smooth function $\Phi : \Omega \to \RR$ is \textbf{radial at x} if there is a function $\Phi_{\mathrm{rad}} : \RR_+ \to \RR$, called \textbf{the radial part} of $\Phi$, such that for any $y \in \Omega$ we have $\Phi(y) =  \Phi_{\mathrm{rad}}(\rho)$ where $\rho = d(x,y)$. Note that the function $\Phi_{\mathrm{rad}}$ is differentiable at $0$ and such that $ \Phi_{\mathrm{rad}}(0) = \Phi(x)$. \\

Let $\Omega, \Phi, \rho, \lambda$ and $x$ as in Lemma \ref{lemmaappendix2}. We denote by $K_x$ the compact group of isometries of $\HH^d$ fixing $x$ and by $\mu_{K_x}(g)$ its Haar measure of unit mass. Because $\HH^d$ is a symmetric space of rank 1, $K_x$ acts transitively on the unit tangent  sphere at $x$. If $y \in \HH^d$ and letting $\rho := d(x,y)$ we thus have
 	$$ \widehat{\Phi}_x(y) := \inte{K_x} \Phi \circ g \ d \mu_{K_x}(g) = \finte{S(x, \rho)}{\Phi}{ d \mu_{S(x,\rho)}}  \ . $$

In particular the function  $\widehat{\Phi}_x$ is radial at $x$ and 
	$$ (\widehat{\Phi}_{x})_{\mathrm{rad}} (\rho) = \finte{S(x, \rho)}{\Phi}{ d \mu_{S(x,\rho)}}  \ . $$

Moreover, since $K_x$ acts by isometries, we also have $\Delta \widehat{\Phi}_x = \lambda \widehat{\Phi}_x$. Therefore, Lemma \ref{lemmaappendix2} reduces to the study of the radial parts of radial at $x$ solution of $\Delta = \lambda \Id$. The following lemma characterises these functions with an ODE. Recall that we denoted by $\mathcal{V}_d(\rho)$ the volume of any ball of radius $\rho$. 
\begin{lemma}
	\label{lemappendix1}
Let $x \in \HH^d$, $\lambda \ge 0$  and $\Omega$ an open set of $\HH^d$ containing a $\rho_0$-neighbourhood of $x$. Let $\Phi$ be a radial at $x$ function satisfying $\Delta \Phi = \lambda \Phi$. Then its radial part solves the following ODE on the interval $]0,\rho[$.
	\begin{equation}
		\label{eqEDOradial}
		 y'' + \frac{\mathcal{V}''_d}{\mathcal{V}'_d} \ y' = - \lambda y \ .
	\end{equation}
\end{lemma}

The proof consists in computing the Laplace operator in polar coordinates and is left to the reader. \\

The second order ODE \eqref{eqEDOradial} has a unique solution smooth at $0$.

\begin{lemma}
	\label{lemmeEDOradial}
	For any $\rho, \lambda >0$ and $y_0 > 0$ there is at most one solution of the ODE \eqref{eqEDOradial} defined on $]0,\rho[$ which can be extended by $y_0$ at $0$ as a $\mathcal{C}^1$ function.
\end{lemma}

\textbf{Proof.} The ODE \eqref{eqEDOradial} being of order 2, we now that it admits at most a space of dimension 2 of solutions defined on $]0, \rho_0[$. Let $y_1$ and $y_2$ be a basis of the vectorial space of the solution of the ODE \eqref{eqEDOradial}. It is classical that the associated Wronskian
$$ W := \det \begin{pmatrix}
	y_1 & y_2 \\ y_1' & y_2'
\end{pmatrix} $$
satisfies the first order ODE
	$$ W' = \frac{\mathcal{V}''_d}{\mathcal{V}'_d} \ W \ . $$ 

But the Wronskian cannot be smooth at $0$ since the function $\frac{\mathcal{V}''_d}{\mathcal{V}'_d}$ is not continuous at $0$. This prevents two independent solutions of \eqref{eqEDOradial} to have $\mathcal{C}^1$ extensions at $0$. \hfill $\blacksquare$ \\

We set $\mathfrak{S}_{\lambda}$ as the unique solution of \eqref{eqEDOradial} defined on $\RR_+^*$ which extends as a $\mathcal{C}^1$ function on $\RR_+$ with $\mathfrak{S}_{\lambda}(0) =1$. \\

We have shown that the following function (which is $\mathcal{C}^1$ at $0$ and values $\Phi(x)$)
	$$ \rho \mapsto \finte{S(x, \rho)}{\Phi}{ d \mu_{S(x,\rho)}} $$
satisfies the ODE \eqref{eqEDOradial} on $]0,\rho[$. By uniqueness, one must have
	$$  \finte{S(x, \rho)}{\Phi}{ d \mu_{S(x,\rho)}} = 	\mathfrak{S}_{\lambda}(\rho) \ \Phi(x) \ , $$ 
which is the desired conclusion. \hfill $\blacksquare$

\bibliographystyle{alpha}

\bibliography{bibliography} 

\end{document}